\documentclass{amsart}

\usepackage{amsthm, amssymb, amsmath, dsfont, mathrsfs, color}
\usepackage[margin=1.2in]{geometry}
\usepackage{graphicx}
\usepackage{subfig}

\newtheorem{theorem}{Theorem}
\newtheorem{corollary}{Corollary}
\newtheorem{lemma}{Lemma}
\newtheorem{proposition}{Proposition}
\theoremstyle{definition}
\newtheorem{definition}{Definition}

\newtheorem{remark}{Remark}

\newcommand*{\E}{\mathbb{E}}

\newcommand*{\N}{\mathbb{N}}
\renewcommand{\P}{\mathbb{P}}

\newcommand*{\R}{\mathbb{R}}

\newcommand{\new}[1]{{\color{red} #1}}

\begin{document}
\title{Asymptotic results of a multiple-entry reinforcement process}	
\author{Caio Alves, Rodrigo Ribeiro and Daniel Valesin }
\date{\today}

\maketitle
\begin{abstract}
We introduce a class of stochastic processes with reinforcement consisting of a sequence of random partitions $\{\mathcal{P}_t\}_{t \ge 1}$, where $\mathcal{P}_t$ is a partition of $\{1,2,\dots, Rt\}$.  At each time~$t$,~$R$ numbers are added to the set being partitioned; of these, a random subset (chosen according to a time-dependent probability distribution) joins existing blocks, and the others each start new blocks on their own.  Those joining existing blocks each choose a block with probability proportional to that block's cardinality, independently. We prove results concerning the asymptotic cardinality of a given block and central limit theorems for associated fluctuations about this asymptotic cardinality: these are proved both for a fixed block and for the maximum among all blocks. We also prove that with probability one, a single block eventually takes and maintains the leadership in cardinality. Depending on the way one sees this partition process, one can translate our results to Balls and Bins processes, Generalized Chinese Restaurant Processes, Generalized Urn models and Preferential attachment random graphs.
\end{abstract}

\section{Introduction}
%!TEX root= ./ms.tex
%In behavioral psychology, reinforcement is a consequence applied that will strengthen an organism's future behavior whenever that behavior is preceded by a specific antecedent stimulus. This strengthening effect may be measured as a higher frequency of behavior (e.g., pulling a lever more frequently), longer duration (e.g., pulling a lever for longer periods of time), greater magnitude (e.g., pulling a lever with greater force), or shorter latency (e.g., pulling a lever more quickly following the antecedent stimulus). There are two types of reinforcement, known as positive reinforcement and negative reinforcement; positive is where by a reward is offered on expression of the wanted behaviour and negative is taking away an undesirable element in the persons environment whenever the desired behaviour is achieved.

The term \textit{reinforcement} in behavioral psychology refers to an effect applied to an organism, triggered by a manifestation of a certain behavior, with the intention of either strengthening or weakening that behavior. In probability, the term describes a large collection of models, such as Balls and Bins models \cite{EP1923, PRR19}, Random Graphs models with Preferential Attachment \cite{ARS17b,BA99} and Reinforced Random Walks \cite{D90}. Although it may be hard to trace its origin, the use of the term reinforcement in the context of random processes has a clear parallel with its meaning in behavioral psychology. 

One of the most important processes in the family of reinforced random processes is the classical P\'olya urn model. In this model, the act of adding balls of a given color with probability proportional to the amount of balls with that color is a reinforcement whose effect is strengthening the presence of that color. In preferential attachment random graphs, the reinforcement strengthens the behavior of choosing high-degree vertices. We do not intend to cover the huge bibliography in the subject,  and refer the interested reader to R. Pemantle's survey on reinforced random processes \cite{P07} and references therein.

In a high level, one of the main goals when dealing with reinforced random processes is to understand the long-run effect of the reinforcement: one wants to check whether the reinforced behavior propagates forever. In the case of models in which the intensity of reinforcement can be tuned via a parameter, one of the main questions is if the model undergoes some sort of phase transition.

In this paper, we introduce and study a class of stochastic processes with reinforcement which can be seen as a meta-model in the sense that by changing the terminology it becomes either an alternate version of Pitman's Chinese restaurant model~\cite{pitman1995exchangeable,pitman2002combinatorial} or a Balls and Bins models or a random (hyper)graph model with preferential attachment. In order to better introduce and discuss our results we will define the model in an abstract setup in Section \ref{sec:model} and then in Section \ref{sec:discussion} we discuss how the model and our results can be interpreted under different perspectives.

\subsection{Definition of the Model}\label{sec:model}

Our model has as parameters a natural number~$R$ and a sequence of probabilities on the $R$-simplex, that is, a sequence of vectors~$(g_t(0),\ldots, g_t(R))_{t \ge 2}$ with~$g_t(r) \ge 0$ for every~$r,t$ and~$\sum_{r=0}^R g_t(r) = 1$ for every~$t$. The dynamics is defined from a sequence~$\mathcal{X}_2,\mathcal{X}_3,\ldots$ of independent random elements of~$\{0,\ldots, R\}$ such that the distribution of each~$\mathcal{X}_t$ is given by the probability mass function~$g_t$. 
The process defines a sequence of random partitions~$\{\mathcal{P}_t\}_{t\in\N}$, where for each~$t$, $\mathcal{P}_t$ is a random partition of $\{1,2,\dots, Rt\}$. We construct these partitions as follows. At time~$t=1$, we start with the set $\{1,\dots, R\}$ and the partition $\mathcal{P}_1$ consisting of a single block equal to~$\{1,\ldots,R\}$, labeled block~1. Assume that the model has been defined up to time~$t$, and that at time~$t$, there are~$N_t$ blocks, labeled from~$1$ to~$N_t$, with~$d_t^{(i)}$ numbers in block~$i$, for~$i \in \{1,\ldots, N_t\}$. Then, at time~$t+1$, from the partition~$\mathcal{P}_t$ of~$\{1,\ldots,Rt\}$ we obtain a partition~$\mathcal{P}_{t+1}$ of~$\{1,\ldots, R(t+1)\}$ as follows:
\begin{itemize}
\item the numbers~$Rt+1,\ldots, Rt+\mathcal{X}_{t+1}$ join blocks that were already present in~$\mathcal{P}_t$, as follows. Independently, each of these numbers joins block~$i\in\{1,\ldots, N_t\}$ with probability
\begin{equation}\label{eq:choose_table}\frac{d_t^{(i)}}{\sum_{j=1}^{N_t} d_t^{(j)}}\end{equation}
(in case~$\mathcal{X}_{t+1} = 0$, this step is skipped); 
\item the remaining numbers~$Rt+\mathcal{X}_{t+1}+1,\ldots, R(t+1)$ are each included in a new block of its own; these new blocks receive labels~$N_t+1,\ldots,N_t+R-\mathcal{X}_{t+1}$ (in case~$\mathcal{X}_{t+1} = R$, this step is skipped).
\end{itemize}
By setting~$d_t^{(i)} = 0$ in case block~$i$ has not yet been started at time~$t$, the above gives rise to a process~$(d_t^{(1)},d_t^{(2)},\ldots)_{t \ge 1}$ on~$(\N_0)^{\N}$ which counts the size of each block in the random partition of the sets $(\{1, \dots, Rt\})_{t \geq 1}$.

%Pitman's original model is recovered  with~$R=1$, and the sequence~$(g_t)$ most often taken as $(g_t(0),g_t(1)) \equiv (p,1-p)$ for some~$p \in (0,1)$ or~$(g_t(0),g_t(1)) = (\frac{1}{t},\frac{t-1}{t})$. Pitman also considered versions of the model in which the probabilities in~\eqref{eq:choose_table} are replaced by other choices, but we do not incorporate this possibility here.

\subsection{Main results} We now present the main results of the paper. In Section \ref{sec:discussion}, we will present a deeper discussion of our results and relate them to other results in the literature. In order to properly state our results we will need some assumptions and definitions.

We make the following assumptions concerning the sequence~$(g_t)_{t\ge 2}$. These assumptions are in force throughout the paper, so we only mention them once here.\\[-.1cm]

\noindent \textbf{Assumptions.} We assume that~$(g_t)_{t \ge 2}$ satisfies
\begin{equation}\label{eq:sum_new}\tag{A.1}
    \sum_{t = 2}^\infty (1-g_t(R)) = \infty.
\end{equation}
By the Borel-Cantelli lemma, this assumption is equivalent to the assumption that infinitely many blocks appear in the process. We also assume that there exists a function $g_\infty: \{0,\ldots, R\} \to [0,1]$ with~$g_\infty(0) < 1$ and 
\begin{equation}\label{eq:assumption}\tag{A.2}
    \sum_{t = 2}^\infty \frac{1}{t}\cdot \sum_{r=0}^R {|g_t(r) - g_\infty(r)|} < \infty,
\end{equation}
which can be seen as a regularity condition for the convergence of $g_t \to g_{\infty}$.

Let~$\beta_t$ be the expected proportion of numbers included in previously existing blocks at time~$t$, that is,
	\begin{equation}\label{eq:def_of_gamma_t}
	\beta_t:= \frac{1}{R}\sum_{r=0}^R r\cdot g_t(r),\qquad t \in \mathbb{N},
	\end{equation}
Also let $\beta$ be its `infinity' counterpart, that is
$$\beta := \frac{1}{R}\sum_{r=0}^R r \cdot g_\infty(r) = \lim_{t \to \infty} \beta_t.$$

We note that~$0 < \beta \le 1$. The case~$\beta =1$, which is equivalent to~$g_\infty(R) = 1$, is somewhat special: in this case, even though it still happens that infinitely many blocks are created (by~\eqref{eq:sum_new}), the expected number of new blocks created at each time step vanishes as~$t \to \infty$. In several of our results, our model undergoes a phase transition at~$\beta = 1$.

Our first result shows that the cardinality of a given block in the random partition grows as~$t^\beta$ multiplied by a random limiting value.
\begin{theorem} [Convergence] \label{thm:convergence}
	For each~$i \in \mathbb{N}$, as~$t \to \infty$,~$d_t^{(i)}/t^\beta$ converges almost surely and in~$L_p$ (for any~$p \in [1,\infty))$ to a strictly positive (and finite) random variable~$\xi^{(i)}$. %Moreover, as~$t \to \infty$,~$\max_i d^{(i)}_t/\phi_t$ converges almost surely and in~$L^p$ (for any~$p \in [1,\infty)$) to a strictly positive (and finite) random variable~$\zeta^*$.
\end{theorem}
Recall that at time $t$ our model gives a random partition of $\{1,2,\dots,Rt\}$. The above theorem reveals a phase transition on the proportion of~$\{1,\ldots, Rt\}$ taken by any fixed set of the random partition. In the regime~$\beta <1$, this proportion vanishes as $t$ goes to infinity, whereas in the regime~$\beta = 1$, it converges to a non-zero limit. 

% Notice that in the regime~$\beta < 1$, the above theorem guarantees that for any~$i \in \N$,
% $$\lim_{t \to \infty} \frac{1}{Rt}\cdot \sum_{j=1}^i d_t^{(j)} = 0 \quad \text{almost surely,}$$
% that is, the density of the first~$i$ sets vanishes as~$t \to \infty$. This is drastically different from the case where~$\beta = 1$, where the density of any given set in the random partition does not vanish; in fact, with some more effort we can also prove:

% \begin{theorem}\label{thm:large_proportion}
% If~$\beta = 1$, then for any~$\varepsilon > 0$ there exists~$i \in \N$ such that, with probability above~$1-\varepsilon$,
% $$1-\varepsilon \le \frac{1}{Rt}\cdot \sum_{j=1}^i d_t^{(j)} \le 1\quad \text{for all }t \ge 1.$$
% \end{theorem}
Our next result provides finer information about the almost sure convergence stated in Theorem~\ref{thm:convergence}. We have a central limit theorem, that is, the random process $\{d^{(i)}_t/t^{\beta}\}_t$ has fluctuations given by a mixed normal distribution around its random limit $\xi^{(i)}$.

\begin{theorem}[Central limit theorem] \label{thm:clt}
	For each~$i \in \mathbb{N}$, 
	$$t^{\beta/2}\cdot \left(\frac{d_t^{(i)}}{t^\beta} - \xi^{(i)} \right) \xrightarrow[\mathrm{(d)}]{t \to \infty} \mu^{(i)},$$
	where~$\mu^{(i)}$ is the distribution of~$W\cdot Z^{(i)}$, where~$W, Z^{(i)}$ are independent random variables,~$W$ is a standard Gaussian and
	$$(Z^{(i)})^2 \sim \begin{cases}\xi^{(i)} &\text{if } \beta < 1;\\ \xi^{(i)}\left(1- \frac{\xi^{(i)}}{R} \right)&\text{if } \beta = 1. \end{cases}$$
\end{theorem}
\iffalse
\begin{remark}
	Recall that when~$\beta = 1$, we have~$\frac{\phi_t}{t} \xrightarrow{t \to \infty} b$ and~$\frac{d_t^{(i)}}{\phi_t} \xrightarrow{t \to \infty} \zeta^{(i)}$ almost surely; together with~$d_t^{(i)} \le Rt$, this gives~$\frac{\zeta^{(i)}b}{R} \in [0,1]$.
\end{remark}
\fi
Since, by Theorem \ref{thm:convergence}, $|d_t^{(i)}/t^{\beta} - \xi^{(i)}|$ converges a.s. to zero as $t$ goes to infinity, Theorem \ref{thm:clt} gives additional information on the rate of convergence. Its statement guarantees that $|d_t^{(i)}/t^{\beta} - \xi^{(i)}|$ rescaled properly has Gaussian fluctuations.

The next result concerns the blocks with the largest cardinality.
\begin{theorem}[Persistent leadership] \label{thm:leader} Almost surely, there exists~$I \in \mathbb{N}$ such that
$$d^{(I)}_t - \max_{i \neq I} d^{(i)}_t \xrightarrow{t \to \infty} \infty.$$
\end{theorem}
The above theorem implies not only that we have a unique block with the maximum cardinality in the random partitions for large enough $t$, but that this block remains far ahead of its ``competition'' forever. This result gives an appealing picture of leadership in the various interpretations of our model. It also has important technical implications. It guarantees in particular that we can extend our central limit theorem to the size of the largest block.
\begin{theorem}[Convergence and central limit theorem for the maximum]
\label{thm:max}
Letting~$I$ be as in Theorem~\ref{thm:leader}, we have
$$\frac{\max_i d_t^{(i)}}{t^\beta} \xrightarrow{t\to \infty} \xi^{(I)} = \sup_{i \in \N} \xi^{(i)} \in (0,\infty),$$
the convergence holding almost surely and in~$L_p$ for any~$p \in [1,\infty)$. 
Moreover,
$$t^{\beta/2}\cdot \left(\frac{\max_i d_t^{(i)}}{t^\beta} - \xi^{(I)} \right) \xrightarrow[\mathrm{(d)}]{t \to \infty} \mu^{*},$$
	where~$\mu^{*}$ is the distribution of~$W \cdot Z^{*}$, where~$W,Z^{*}$ are independent,~$W$ is a standard Gaussian and~$$(Z^{*})^2 \sim  \begin{cases}\xi^{(I)}&\text{if } \beta < 1;\\ \xi^{(I)}\left(1- \frac{\xi^{(I)}}{R}\right) &\text{if }\beta = 1.\end{cases}$$
\end{theorem}

\subsection{Discussion of the results and related work}\label{sec:discussion}
In this section we discuss our results in more detail and show how they can be interpreted under the terminology of other classical models. \textit{En passant}, we mention how our model can be seen as a generalization of other well-known reinforced random processes.

\subsubsection*{\textbf{Cardinality distribution}} In this paper, we are concerned with cardinalities of specific blocks, as well as the maximum cardinality. A natural follow-up question would then be: how are the cardinalities, taken as a whole, distributed? More specifically, what is the proportion of blocks with a given cardinality~$k$? Calculations similar to the proof of Theorem 8.3 of \cite{RemcoBook16} show that, when~$\beta < 1$, this cardinality distribution has power-law with exponent~$1 + 1/\beta \in (2, \infty)$. When~$\beta = 1$ however, the speed of convergence of~$\beta_t$ to~$\beta$ directly influences the cardinality distribution, some sequences might not necessarily yield processes with power-law cardinality distribution. We plan to study these properties in depth in the future.

\subsubsection*{\textbf{Universality}} Note that our limit random variables $\xi^{(i)}$ depend on the choice of the sequence of vectors $(g_t)_{t\ge 2}$, where $g_t$ is a point of the $R$-dimensional simplex. On the other hand, qualitatively, our results depend on a functional of the limiting point $g_{\infty}$. Moreover, this functional has a simple formula: it is just the expected value associated to $g_{\infty}$ divided by $R$. In other words, there exists universality on the asymptotic behavior of the observables $d_t^{(i)}$'s as well as their maximum: the order of magnitude of those observables is determined by a functional of the limiting point $g_{\infty}$. Thus, any sequence of points $(g_t)_{t\ge 2}$ in the $R$-dimensional simplex satisfying assumptions (\ref{eq:sum_new}) and (\ref{eq:assumption}) produces the same asymptotic behavior, though the distributions of limiting random variables $\xi^{(i)}$'s may depend on the choice of sequence $(g_t)_{t\ge 2}$.

\subsubsection*{\textbf{Balls and Bins models}} In this case there is a natural correspondence between the terminology of random partitions and a balls and bins model with the feature that at each step a random number of new bins is added by the process. The bins are the blocks in the random partition and the natural numbers corresponds to labels to the balls:  interpret each natural number as a ball labelled with that number, and a block of a partition as a bin. Hence, the state of the process at time~$t$ describes an assignment of~$Rt$ balls (labelled from~$1$ to~$Rt$) into bins. At each time step,~$R$ new balls and a random number of bins join the process.

In the balls and bins scheme, one of the central questions regards \textit{dominance} and \textit{monopoly}. Dominance means that after some (random) time there will be a single bin with more balls than the others, whereas monopoly means that after some (random) time only one bin will receive all the balls. In these terms, Theorem \ref{thm:leader} states that we do observe dominance. However, Theorem \ref{thm:convergence} guarantees there is no monopoly, since it guarantees that the number of balls in any bin goes to infinity almost surely.

Still in the context of balls and bins models, one feature of our results stands out: the dependence of the quantity $\beta$. Our results state that under assumptions~\eqref{eq:sum_new} and \eqref{eq:assumption}, the growth rate of the number of balls inside each bin is determined by $\beta$, which in this context reads as the limiting average proportion of balls added to already existing bins.

\subsubsection*{\textbf{Urn model: One vs. All}} If one desires to keep track of the size of a single block, then our model becomes the usual urn model with one urn and two colors, red and green, and the additional feature of \textit{immigration}.
The model consists in a urn which at time~$t$ contains a collection of~$Rt$ balls colored either red or green: red balls represent numbers included in a chosen block of the partition, and green balls represent every other number. At each time $t$, $\mathcal{X}_t$ new balls are put into the urn with colors chosen in the usual fashion without reposition, together with $R-\mathcal{X}_t$ green balls.

In this context, $\beta_t$ becomes the average proportion of non-immigrant balls added at time $t$. Our results then make the effect of immigration clear: when the asymptotic proportion of immigrant balls being added is zero (in which case $\beta = 1$), the limiting proportions of red and green balls are comparable in the sense that
$$
\lim_{t\to\infty}\frac{\# \text{red balls at time }t}{\# \text{green balls at time }t} = \frac{\xi^{(i)}}{R-\xi^{(i)}},
$$
which is a strictly positive random variable. On the other hand, when the asymptotic proportion of immigrant balls being added is positive, $\beta < 1$,  then
$$
\lim_{t\to\infty}\frac{\# \text{red balls at time }t}{\# \text{green balls at time }t} = 0.
$$ 
We refer the interested reader in the subject of urn models with immigration to \cite{PRR19} of E. Pek{\"o}z, A. R{\"o}llin and N. Ross about urn model with random immigration and references therein.

\subsubsection*{\textbf{Generalized Chinese Restaurant Process}}Given~$\alpha \leq 1$ and~$\theta \geq -\alpha$, the Generalized Chinese Restaurant Process, denoted by $\mathrm{GCRP}(\alpha, \theta)$, is a model introduced by J. Pitman \cite{pitman1995exchangeable} as a generalization of Ewens formula for sampling \cite{ewens1972}. In this process, a configuration at a given time~$t$ is given by a finite collection of tables, each table having a finite number of customers and a natural number as index. The process is constructed inductively, at time $t+1$ a customer arrives at the restaurant and either: sits in the already existing $i$-th table with probability 
$$
\frac{\# (\text{number of customers in the } i \text{-th table at time }t) - \alpha}{t+ \theta}
$$
or sits in a new table with probability 
$$
\frac{\alpha \cdot (\# \text{number of tables at time }t) + \theta}{t+ \theta}.
$$
When $\alpha = 0$ and $\theta = 1$, this is the usual Chinese Restaurant Process.

If we let $T_t$ be the number of tables at time $t$, $\mathrm{GCRP}(\alpha, \theta)$ presents different behaviors which can be summarized as follows:
\begin{enumerate}
	\item \textit{Finite number of tables:} Choosing $\alpha < 0$ and $\theta = -m\alpha$, for some $m$, then $T_t \to m$ and after a finite random time the model behaves like a balls and bins model with $m$ bins;
	\item \textit{Logarithmic growth:} Setting $\alpha=0$ and $\theta > 0$, we have $T_t$ of order $\log t$;
	\item \textit{Polynomial growth:} For $\alpha >0$ and $\theta > -\alpha$, $T_t$ grows like $t^{\alpha}$. 
\end{enumerate}
The interested reader may see \cite{pitman2002combinatorial} for an overview of such regimes.

Under the terminology of the $\mathrm{GCRP}$, our model becomes a quite natural alternate version of it and reads as follows: at each time  $t$, $R$ customers arrive to the Chinese restaurant; $\mathcal{X}_t$ of them sit in already set tables, choosing their tables with probability proportional to the number of customers at each table; and $R-\mathcal{X}_t$ customers sit each one on a new table. %Notice that, choosing $R=1$ and $g_t$ as~$g_t(1) = 1-g_t(0) = \frac{1}{t+1}$ we obtain the Chinese Restaurant Process.

Regarding the growth of the number of tables, by choosing the sequence $(g_t)_{t \ge 2}$ properly all the behaviors described above may be achieved. On the other hand, in contrast to $\mathrm{GCRP}$ with $\alpha >0$, in our case the growth of the number of tables is driven by the process $\{\mathcal{X}_t\}_{t \ge 2}$ and not only $\{T_t\}_{t \ge 1}$. This distinguishes our model from the $\mathrm{GCRP}$ even when both have polynomial growth in the number of tables, and puts it into a new category of Chinese restaurant processes.

When dealing with $\mathrm{GCRP}$ the central questions regard the number of tables of a given size (the amount of customers sitting at the table) and the size of the tables, see [favaro2015,favaro2018,OPR20]. In this context our results regard the size of the tables  (Theorem \ref{thm:convergence}), the size of the largest table (Theorem \ref{thm:max}), the existence of a single largest table (Theorem \ref{thm:leader}). All of these results are sensitive to the asymptotic proportion of new tables being added, $1-\beta$.

\subsubsection*{\textbf{Random Graph model}} In the context of random graph models whose vertex set may increase with time, our model can be seen as a building block and give rise to different random graph models. Here each block of the random partition $\mathcal{P}_t$ can be interpreted as a vertex. The new blocks which may be introduced at each step are new vertices. Then, in order to generate a graph from this configuration, one has to impose a rule of connection which is a rule to link the blocks in the random partition. 

In order to better illustrate this \textit{rule of connection}, let us define a rule which introduces geometry to our model. In this analogy, the starting partition $\mathcal{P}_1 = \{\{1,2,\dots, R\}\}$ will correspond to a starting graph consisting of one vertex with $Rd$ loops. Assume then that the graph generated from the partition~$\mathcal{P}_t$ is already constructed, and let us inductively define the graph associated to time~$t + 1$. Let $H_d$ be a $d$-regular graph with $R$ vertices labeled $h_1,h_2,\dots, h_R$ (notice that the existence of such $H_d$ implies $Rd$ is even). We also denote by $B_{Rt+i}$, with $i \in \{1,2,\dots, R\}$, the block of the partition~$\mathcal{P}_{t + 1}$ which receives the number $Rt+i$ at time~$t + 1$. Then our rule of connection works as follows: we add an edge between blocks $B_{Rt+i}$ and $B_{Rt+j}$ if and only if $h_i$ and $h_j$ are connected in $H_d$. When $B_{Rt+i} = B_{Rt+j}$ we add a loop to the corresponding vertex. Note that, since the graph~$H_d$ associated to the rule of connection is $d$-regular, the degree of a vertex corresponds to $d$ times the cardinality of its associated block. Therefore, we can translate our results to the degree of a given vertex and the maximum degree in this graph process. Given the interest in graph observables related to counting specific subgraphs such as triangles, see \cite{JM15,KH02, ORS18}, choosing $H_d$ and $g_t$'s properly may lead to interesting models.

\begin{figure}[!tbp]
  \centering
  \includegraphics[width=0.7\textwidth]{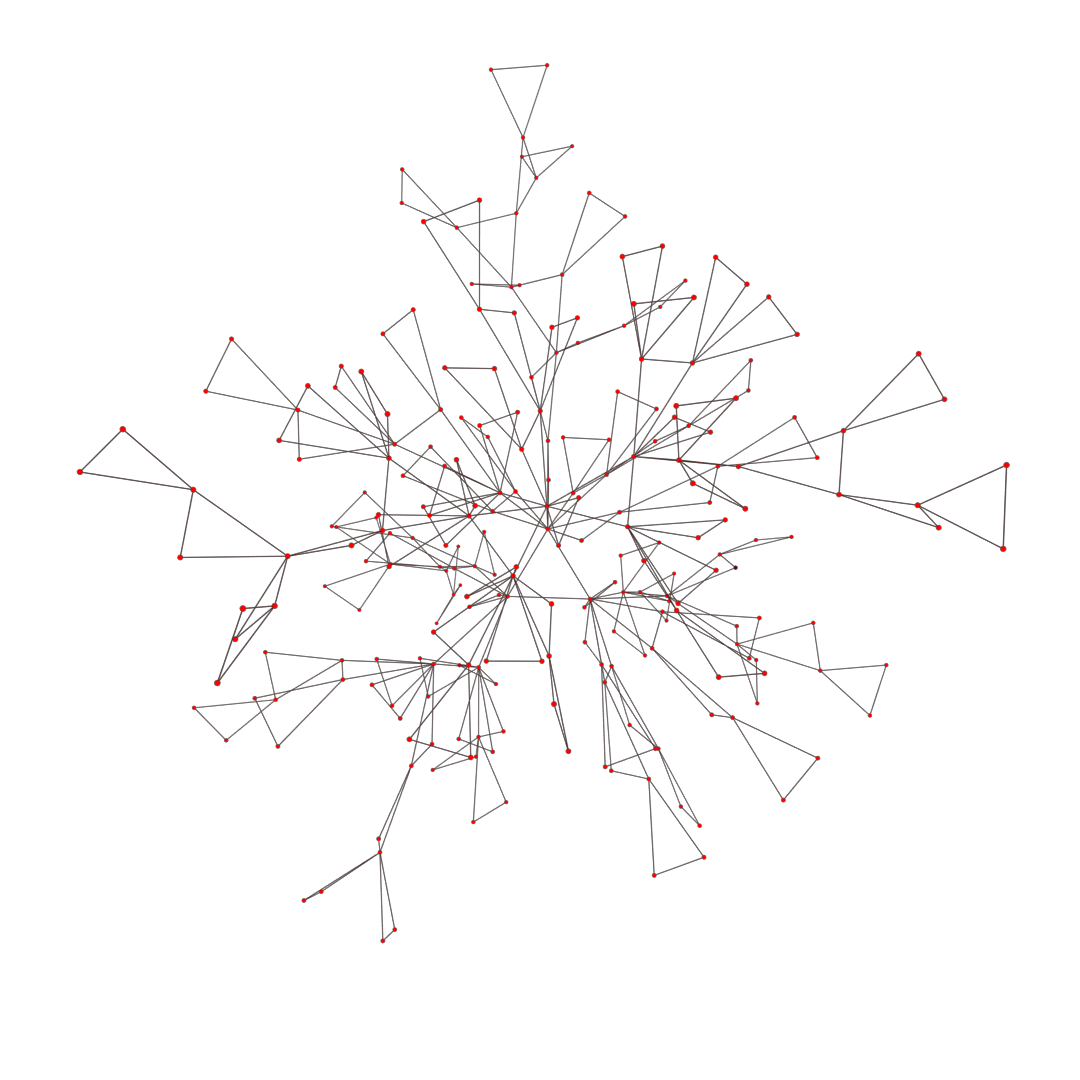}\label{fig:K_3_graph}
  \caption{A simulation of the resulting graph with circa $250$ vertices, when~$H_d$ is the triangle~$K_3$ and~$g_t(1) \equiv 1$. In this case $\beta = 1/3$.}
\end{figure}

For instance, by taking $R=2$, $d=1$ and $H_d$ as a single edge between the two vertices we can recover preferential attachment graphs. By setting $g_t(1) = 1$ this is an instance of the  A. Barab\'asi and R. \'Albert random tree \cite{BA99}. Whereas, by taking $g_t(1) = 1- g_t(2) = f(t)$ our model becomes the preferential attachment random graph recently introduced in \cite{ARS17b}, in which one can add edges between already existing vertices.

By choosing a different rule of connection, our model can generate random hypergraphs as well. In this variant of the model, the vertices will still be the blocks of the partition at time~$t$, but now we add a single hyper-edge between all the vertices (or blocks) which have increased its degree at each time-step, instead of connecting them based on the chosen graph~$H_d$. This random hypergraph model behaves similarly to the Hollywood model introduced by H. Crane and W. Dempsey in \cite{crane2017} 

\begin{figure}[!tbp]
  \centering
  \includegraphics[width=0.7\textwidth]{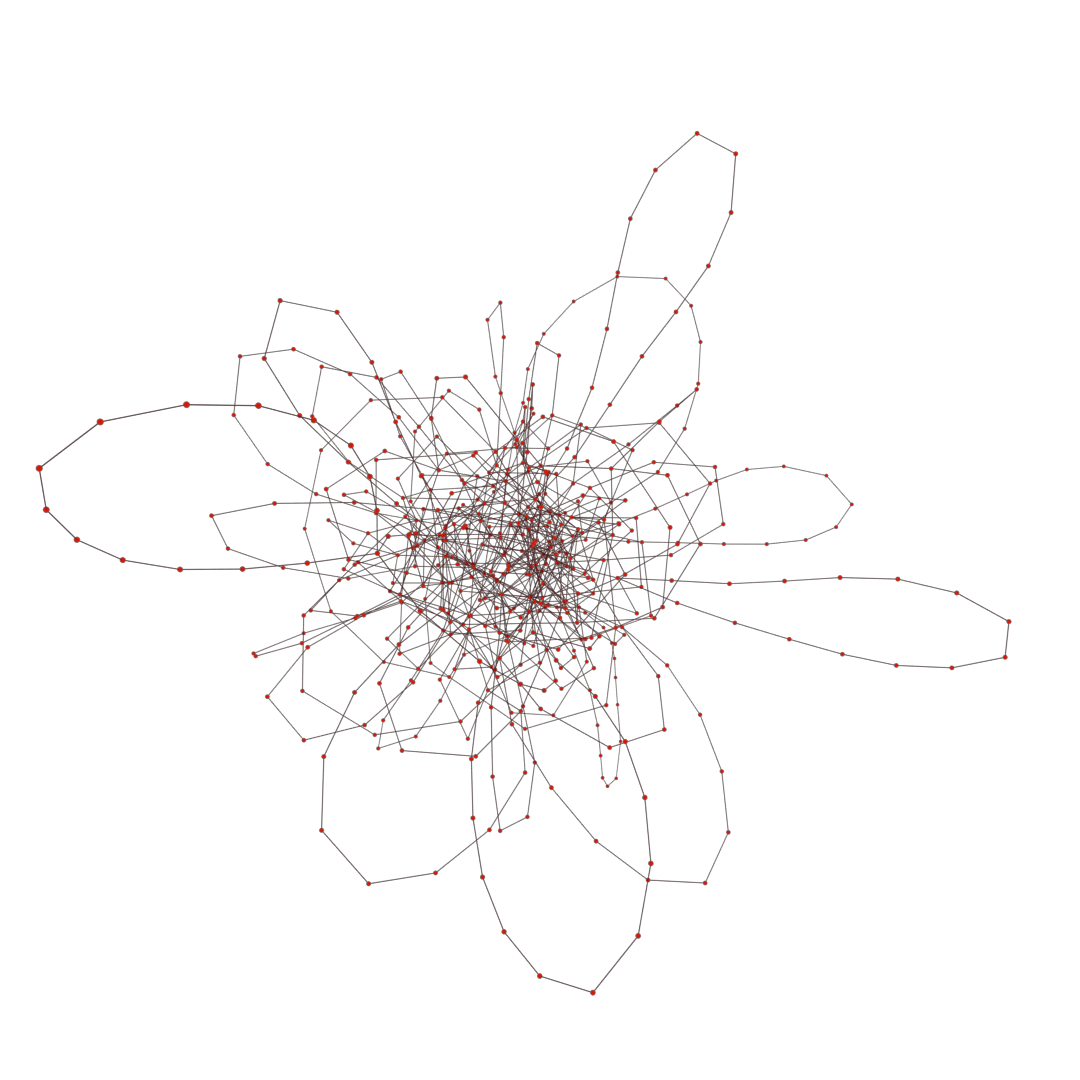}\label{fig:C_15_graph}
  \caption{A simulation of the resulting graph with circa $250$ vertices, when~$H_d$ is the cycle~$C_{15}$, $g_t(1) \equiv 0.7$, and~$g_t(15) \equiv 0.3$. In this case $\beta = 26/75$.}
\end{figure}

Under these terminologies, Theorem \ref{thm:leader} guarantees that, for large~$t$, there exists only one vertex with maximum degree, which keeps its leadership forever. Moreover, Theorem \ref{thm:max} guarantees a phase transition on the order of the maximum degree. In the regime $\beta < 1$ we have that 
$$
\lim_{t \to \infty}\frac{\text{maximum degree at time }t}{t} = 0, \text{ a.s.}
$$
On the other hand, for $\beta = 1$, the model achieves a maximum degree of linear order, that is
$$
\lim_{t \to \infty}\frac{\text{maximum degree at time }t}{t} > 0, \text{ a.s.}
$$

\subsection{Organization of the paper}
In Section~\ref{sec:conv}, we will study a process given by the cardinality of a given block divided by a normalizing function of~$t$. This will allow us to use martingale arguments and results: we use Freedman's inequality to obtain an exponential upper bound for the tail of the normalized cardinality; M\'ori's martingale~\cite{M05} to study moment bounds and $L_p$-convergence of the normalized cardinality; and a generating function argument based on~\cite{S18} to prove a.s. positivity of the limiting random variable, thus proving Thorem~\ref{thm:convergence}.

In Section~\ref{sec:clt}, we will prove the CLT result of Theorem~\ref{thm:clt}, based on an application of Corollary~3.5 of \cite{MartingaleCLTbook2}. In order to show that our process satisfies the hypotheses, we will study the convergence of sums of squared increments of the normalized process.

In Section~\ref{sec:leader} we will prove the persistent leadership result, Theorem~\ref{thm:leader}, based on a Lyapunov function argument stated in \cite{menshikov2008urn}. Finally, in Section~\ref{sec:maximum}, we show that the maximum cardinality also satisfies results similar to the ones proven for the cardinality of a given block, thus proving Theorem~\ref{thm:max}. Here, it will be crucial to use Theorem~\ref{thm:leader} in order to show that, for large~$t$, the maximum cardinality behaves like the cardinality of a given vertex.

\subsection{Notation and conventions}
We write~$\N = \{1,2,\ldots\}$ and~$\N_0 = \{0,1,\ldots\}$.

Keeping consistency with the notation~``$d^{(i)}_t$'', we use super-indices (inside parentheses) to denote table indices (as well as other parameters, as in Definition~\ref{def:phi_t_m} below) and we use sub-indices for time or time-related indices. In general, for an expression with a sub-index, such as~$x_t$, we write
$$\Delta x_t := x_{t+1} - x_t.$$

Given~$n,m \in \mathbb{N}_0$ and~$p \in [0,1]$, we denote by~$\mathscr{B}(n,p,m)$ the probability that a Binomial($n,p$) random variable is equal to~$m$.

\section{Concentration and convergence}\label{sec:conv}
The goal of this section is to prove Theorem \ref{thm:convergence} together with exponential tail bounds for the normalized cardinality of a given block $i$. In order to do that, we need to prove several intermediate results. The almost sure convergence of $\{d_t^{(i)}/t^{\beta}\}_{t\ge 1}$ will follow from martingale theory. Then, convergence in the $L_p$ sense and positiveness of the limiting random variable $\xi^{(i)}$ will require more work.
\subsection{Martingales and concentration}
In this part we will make sure $\{d_t^{(i)}/t^{\beta}\}_{t\ge 1}$ converges almost surely by finding a related martingale. In order to do this the following lemma is vital.
\begin{lemma}\label{lem:expectation_delta_d}
For each~$i, t \in \mathbb{N}$ the following identity holds
\begin{equation}\label{eq:basic_exp_delta}
\E[\Delta d^{(i)}_t \mid \mathcal{F}_t] = \beta_{t+1} \cdot \frac{d_t^{(i)}}{t}, \text{almost surely.}
\end{equation}
\end{lemma}
\begin{proof}
The proof follows by direct computation
\begin{align*}
\E[\Delta d_t^{(i)}\mid \mathcal{F}_t] &= \sum_{r=0}^R r \cdot \P(\Delta d^{(i)}_t = r \mid \mathcal{F}_t) \\&= \sum_{r=0}^R r\cdot \sum_{u=0}^R g_{t+1}(u) \cdot \mathscr{B}\left(u,\frac{d_t^{(i)}}{Rt},r\right) = \sum_{u=0}^R g_{t+1}(u)\cdot u \cdot \frac{d_t^{(i)}}{Rt} = \beta_{t+1}\cdot \frac{d_t^{(i)}}{t}.
\end{align*}
\end{proof}
Recall by the definition of the process that the time the $i$-th block is introduced is random and depends on the sequence $\{\mathcal{X}_t\}_{t\ge 2}$, for this reason we need the following definitions: for each~$i \in \N$, define the stopping time
\begin{equation*}
\tau^{(i)} := \inf\{t \ge 1:\;d^{(i)}_t > 0\}.
\end{equation*}
Also define, for~$i,n \in \N$, with $n\ge i$, the conditional probability measures
\begin{equation*}
\P_n^{(i)}(\;\cdot\;) := \P(\;\cdot \; \mid \tau^{(i)} = n).
\end{equation*}

\begin{remark}
In the above definition and in the rest of the paper, whenever we fixed pair~$i,n$ and refer to the event~$\{\tau^{(i)} = n\}$, we assume tacitly that~$i,n$ are such that this event has positive probability. Note in particular that~$i=1$ then also forces~$n = 1$.
\end{remark}
The next result ensures the $i$-th block is eventually created, for all $i\in \mathbb{N}$. In other words, infinitely many blocks are created.
\begin{lemma}
Almost surely,~$\tau^{(i)} < \infty$ for all~$i \in \N$.
\end{lemma}
\begin{proof}
Notice that the probability that at least one block is created at time~$t$ is~$1-g_t(R)$, so by the assumption~\eqref{eq:sum_new} and the Borel-Cantelli Lemma, almost surely the number of blocks created is infinite.
\end{proof}
In the next results we will show that the cardinality of a given block properly normalized is a (sub)martingale. The factor of normalization at each time together will proper notation is introduced below.
\begin{definition}
	Define the sequence~$(\phi_{t})_{t \ge 1}$ by
	$$\phi_1 := 1 \quad \text{ and }\quad \phi_t := \prod_{s=2}^{t} \left(1+ \frac{\beta_{s}}{s-1} \right),\; t \ge 2 $$
and for each~$i \in \N$, define the process
\begin{equation*}
X_t^{(i)} := \frac{d_t^{(i)}}{\phi_t},\qquad t \ge 1.
\end{equation*}
Note that~$X^{(i)}_t = 0$ for~$t < \tau^{(i)}$.
\end{definition}
The next lemma provides the right order of magnitude of the sequence $(\phi_t)_{t\ge 1}$.
\begin{lemma} \label{lem:convergence_of_phi}
There exists~$b > 0$ such that
\begin{equation}\label{eq:gamma_and_b}
\lim_{t \to \infty} \frac{\phi_t}{t^\beta} = b.
\end{equation}
\end{lemma}
\begin{proof}
Letting~$\mathcal{E}(x):= \log(1+x)-x$ for~$x > -1$, we have, for~$t \ge 2$,
\begin{align*}
\phi_t &= \exp\left\{\sum_{s=2}^t \left(\frac{\beta_s}{s-1}+ \mathcal{E}\left( \frac{\beta_s}{s-1}\right)\right) \right\}\\
&=\exp\left\{\beta \cdot  \log(t) + \beta\cdot \left(\sum_{s=2}^{t} \frac{1}{s-1} - \log(t)\right) + \sum_{s=2}^t \frac{\beta_s - \beta}{s-1} + \sum_{s=2}^t \mathcal{E}\left(\frac{\beta_s}{s-1}\right) \right\},
\end{align*}
then
$$\frac{\phi_t}{t^\beta} = \exp\left\{\beta\cdot \left(\sum_{s=2}^{t} \frac{1}{s-1} - \log(t)\right) + \sum_{s=2}^t \frac{\beta_s - \beta}{s-1} + \sum_{s=2}^t \mathcal{E}\left(\frac{\beta_s}{s-1}\right) \right\}. $$
Now, it is well known that~$\sum_{s=1}^{t-1} \frac{1}{s} - \log(t) $
converges as~$t \to \infty$ to the Euler-Mascheroni constant. Next, the series~$\sum_{s=2}^\infty \frac{\beta_s-\beta}{s-1}$ is convergent due to~\eqref{eq:assumption}. Finally, the series~$\sum_{s=2}^\infty \mathcal{E}\left(\frac{\beta_s}{s-1}\right)$ is convergent since~$\mathcal{E}(x) = o(x^2)$ when~$x\to 0$. 
\end{proof}
As a consequence of Lemma \ref{lem:expectation_delta_d}, we will show that the normalized cardinality of a given block is a (sub)martingale.
\begin{lemma} \label{lem:x_mart}
For any~$i \in \N$, the process~$\{X^{(i)}_t\}_{t \ge 1}$ is a submartingale under~$\P$, and the process~$\{X^{(i)}_t\}_{t \ge n}$ is a martingale under~$\P_n^{(i)}$.
\end{lemma}
\begin{proof}
For the first statement, we start noting that, for any~$i,t$,
\begin{align*}
\E[X_{t+1}^{(i)}\mid \mathcal{F}_t] &= \mathds{1}\{\tau^{(i)} \le t\}\cdot \frac{d_t^{(i)}+\E[\Delta d_t^{(i)}\mid \mathcal{F}_t]}{\phi_{t+1}} + \mathds{1}\{\tau^{(i)}>t\}\cdot \frac{\P(\tau^{(i)} = t+1 \mid \mathcal{F}_t)}{\phi_{t+1}},
\end{align*}
where the expression for the second term on the right-hand side follows from the fact that~$d^{(i)}_{t+1} = 1$ if~$\tau^{(i)} = t+1$. Now, we clearly have
$$\mathds{1}\{\tau^{(i)}>t\}\cdot \frac{\P(\tau^{(i)} = t+1 \mid \mathcal{F}_t)}{\phi_{t+1}} \ge 0 = \mathds{1}\{\tau^{(i)}>t\}\cdot X_t^{(i)},$$
and, using~\eqref{eq:basic_exp_delta}, we also have
$$\mathds{1}\{\tau^{(i)} \le t\}\cdot \frac{d_t^{(i)}+\E[\Delta d_t^{(i)}\mid \mathcal{F}_t]}{\phi_{t+1}} = \mathds{1}\{\tau^{(i)} \le t\}\cdot \frac{d_t^{(i)}+\beta_{t+1}\cdot d_t^{(i)}/t}{\phi_{t+1}} = \mathds{1}\{\tau^{(i)} \le t\}\cdot X_t^{(i)}. $$
This completes the proof of the first statement. For the second statement we write, note that for~$t \ge n$, we have,~$\P^{(i)}_n$-almost surely,
$$\E^{(i)}_n[\Delta d_t^{(i)}\mid \mathcal{F}_t] = \E[\Delta d_t^{(i)}\mid \mathcal{F}_t] \stackrel{\eqref{eq:basic_exp_delta}}{=}   \frac{\beta_{t+1}\cdot d_t^{(i)}}{t},$$
thus
$$\E^{(i)}_n[X_{t+1}^{(i)}\mid \mathcal{F}_t] = \frac{d_t^{(i)}+ \beta_{t+1}\cdot d_t^{(i)}/t}{\phi_{t+1}} = X_t^{(i)},$$
which concludes the proof.
\end{proof}

For the proof of the coming proof and several other points in the paper, it will be useful to write
\begin{equation}\label{eq:a_and_b}
\Delta X^{(i)}_t = A_t - B_t
\end{equation}
where
\begin{equation*}
A_t:= \frac{\Delta d^{(i)}_t}{\phi_{t+1}},\qquad B_t:= X^{(i)}_t \cdot \frac{\Delta \phi_t}{\phi_{t+1}} = X^{(i)}_t \cdot \left(\frac{\phi_{t+1}}{\phi_t} - 1\right) = X^{(i)}_t \cdot \frac{\beta_{t+1}}{t}.
\end{equation*}
It will be useful to control how large the cardinality of a given block can be. In the next lemma we prove exponential tail bounds for $d_t^{(i)}/\phi_t$. The result is the first consequence of the fact that the process $\{X_t^{(i)}\}_{t\ge 1}$ is a martingale under $\P_n^{(i)}$.
\begin{lemma} \label{l:upperbounddeg}
There exist~$c_1,C_1 > 0$ such that, for any~$i,n \in \N$ and any~$\alpha > 0$, 
\begin{equation}
\label{eq:upperbounddeg}
\P^{(i)}_n\left(       \exists t \ge n:\; d_t^{(i)} \ge \alpha\cdot \frac{\phi_t}{\phi_n}     \right) \leq C_1\exp\{    -c_1 \alpha   \}.
\end{equation}
\end{lemma}
\begin{proof}
Fix~$i,n \in \N$. We start with some bounds involving~$\Delta X^{(i)}_t$, for~$t \ge n$. First, using~\eqref{eq:a_and_b} and the facts that~$\phi_{t+1} \ge \phi_t$,~$\Delta d_t^{(i)} \le R$,~$\beta_{t+1}\le 1$ and~$d_t^{(i)} \le Rt$, we have
\begin{equation}\label{eq:bound_triangle_z}
|\Delta X^{(i)}_t| \le \frac{R}{\phi_{t+1}} + \frac{d_t^{(i)}\cdot \beta_{t+1}}{\phi_t\cdot t} \le \frac{2R}{\phi_t}.
\end{equation}
Next, again using~\eqref{eq:a_and_b} and the inequality~$(a+b)^2 \le 2a^2 +2b^2$, we have
\begin{align*}
\E^{(i)}_n[(\Delta X_t^{(i)})^2 \mid \mathcal{F}_t] &\le 2\cdot \frac{\E^{(i)}_n[(\Delta d_t^{(i)})^2 \mid \mathcal{F}_t]}{(\phi_{t+1})^2} + 2\cdot \left(\frac{d_t^{(i)}}{t\cdot \phi_t}\right)^2\\
&\le \frac{2R^2}{(\phi_{t+1})^2}\cdot \P_n^{(i)}(\Delta d_t^{(i)}\neq 0\mid \mathcal{F}_t)+ 2R\cdot \frac{d_t^{(i)}}{t\cdot (\phi_t)^2}.
\end{align*}
Using the bound~$\P^{(i)}_n(\Delta d^{(i)}_t \neq 0\mid \mathcal{F}_t) \le R \cdot \frac{d^{(i)}_t}{Rt} = \frac{d^{(i)}_t}{t}$ and again using~$\phi_{t+1} > \phi_t$, we obtain
\begin{equation}\label{eq:bound2_triangle_z}\E^{(i)}_n[(\Delta X_t^{(i)})^2 \mid \mathcal{F}_t] \le 4R^2 \cdot \frac{d_t^{(i)}}{t\cdot (\phi_t)^2}.\end{equation}

We now fix~$\alpha$. It is sufficient to prove the inequality~\eqref{eq:upperbounddeg} for~$\alpha$ large enough (by increasing~$C_1$ if necessary), so we assume that~$\alpha > R$. Define 
$$\eta:= \inf\left\{t \ge n: \;X^{(i)}_t \ge \frac{\alpha}{\phi_n}\right\}$$
and
$$Z_t := X^{(i)}_{t \wedge \eta} - X^{(i)}_n,\quad t \ge n,$$
which is a martingale under~$\P^{(i)}_n$ by Lemma~\ref{lem:x_mart}. It is worth noting at this point that the case~$i = 1$ is somewhat special: it implies that~$n = 1$ also, and~$X_1^{(1)} = R/\phi_1 = R$. On the other hand, if~$i > 1$ we have~$X^{(i)}_n = 1/\phi_n$. We thus obtain
\begin{equation}\label{eq:z_more}
Z_t \ge X^{(i)}_{t \wedge \eta} - \frac{R}{\phi_n},\quad t \ge n.\end{equation}
Next, from~\eqref{eq:bound_triangle_z} it follows that
\begin{equation}\label{eq:other_bound_z}|\Delta Z_t|\le \frac{2R}{\phi_t} \le \frac{2R}{\phi_n},\quad t \ge n. \end{equation}
Finally, from~\eqref{eq:bound2_triangle_z} it follows that
\begin{equation*}
 \sum_{t \ge n}\E_n^{(i)}[(\Delta Z_t)^2 \mid \mathcal{F}_t] \le \sum_{t \ge n} \mathds{1}\{ \eta \ge t\}\cdot 4R^2\cdot \frac{d_t^{(i)}}{t\cdot (\phi_t)^2} \le 4R^2\alpha \cdot \sum_{t \ge n} \frac{1}{t \cdot \phi_t},
\end{equation*}
since~$d_t^{(i)}/\phi_t = X_t^{(i)} \le \alpha$ on~$\{\eta > t\}$. Then, using the fact that~$\phi_t/t^\beta \to b$ we obtain
\begin{equation}\label{eq:other_bound2_z}
\sum_{t \ge n}\E_n^{(i)}[(\Delta Z_t)^2 \mid \mathcal{F}_t] \le C\alpha\cdot \frac{1}{\phi_n}
\end{equation}
for some~$C > 0$ that does not depend on~$i$,~$n$ or~$\alpha$. Equations~\eqref{eq:other_bound_z} and~\eqref{eq:other_bound2_z} allow us to apply Freedman's inequality, Theorem~\ref{thm:freedman} in the Appendix, with~$\lambda = (\alpha-R)/\phi_n$,~$K:= 2R/\phi_n$ and~$\sigma^2:= C\alpha/\phi_n$, to bound:
\begin{align*}
\P_n^{(i)}\left(\exists t \ge n:\;d_t^{(i)} \ge \alpha \cdot \frac{\phi_t}{\phi_n} \right) &\stackrel{\eqref{eq:z_more}}{\le} \P^{(i)}_n \left(\exists t \ge n: Z_t \ge \frac{\alpha - R}{\phi_n}\right)\\
& \le \exp \left\{-\frac{\left(\frac{\alpha-R}{\phi_n}\right)^2}{2\cdot \frac{C\alpha}{\phi_n}+ \frac23\cdot \frac{2R}{\phi_n} \cdot \frac{\alpha-R}{\phi_n} }\right\} \le \exp\left\{-\frac{(\alpha - R)^2}{2C\alpha\phi_n + \frac43R  (\alpha- R)} \right\}.
\end{align*}
It is now not hard to see (using~$\phi_n \ge 1$) that there exist constants~$C_1,c_1 > 0$ (not depending on~$i,n,\alpha$) such that the right-hand side above is smaller than~$C_1 \exp\{-c_1\alpha/\phi_n\}$.
\end{proof}

\subsection{Higher moments and positivity of the limit} By Lemma \ref{lem:x_mart}, the cardinality of a given block, when properly normalized and considered under the right measure, is a martingale. In this case, almost sure convergence comes from the fact that we have a positive martingale. This however does not imply a crucial property present in Theorem~\ref{thm:convergence}, namely, the fact that the limiting random variable is strictly positive. 

Under our general setting, this is another point in which our model imposes technical difficulties. In other classical and related models, see \cite{M05}, the positivity of the limiting random variable comes naturally from the fact that such limit has a known distribution (a beta distribution in many contexts). Moreover, in many contexts, see \cite{M05}, one of the reasons why it is possible to obtain a explicit formula for the limiting random variable is that at each step only one number is added to the set to be partitioned. Thus, one of the main goals of this section is to derive positivity of the limiting random variable $\xi^{(i)}$.

Our approach was inspired by \cite{M05} and part of it relies on the so-called Mori's martingales, whose definition we give below.

\begin{definition} \label{def:phi_t_m}
For each~$m \in \N$ with~$m \ge 2$, define the sequence~$(\phi^{(m)}_t)_{t \ge 1}$ by
$$\phi_1^{(m)} = 1\qquad \text{and}\qquad\phi_{t}^{(m)} := \prod_{s=2}^{t} \left(1+ m\cdot \frac{\beta_{s}}{s-1} +\frac{R(4m)^R}{(s-1)^2}\right),\; t\ge 2 $$
and for each~$i \in \N$, define the process
$$Y^{(i,m)}_t:= \frac{1}{\phi_t^{(m)}}\cdot {d_t^{(i)}+m-1 \choose m}, \qquad t \ge 1, $$
with the convention that~$\binom{a}{b} = 0$ when~$a < b$ (so that~$Y^{(i,m)}_t > 0$ if and only if~$t \ge \tau^{(i)}$).
\end{definition}
We then have
\begin{lemma} \label{lem:Y_super}
For any~$i,n,m$, we have that under~$\mathbb{P}^{(i)}_n$, the process~$\{Y_t^{(i,m)}\}_{t \ge n}$ is a positive supermartingale.
\end{lemma}
\begin{proof}
We start noting that, on the event~$\{\tau^{(i)}=n\}$, for~$t \ge n$,
\begin{align*}
{d_t^{(i)} + \Delta d_t^{(i)}+m-1 \choose m} &= \sum_{r=0}^R \mathds{1}\{\Delta d_t^{(i)} = r\}\cdot {d_t^{(i)} + r+m-1 \choose m}\\
&= \sum_{r=0}^R \mathds{1}\{\Delta d_t^{(i)} = r\}\cdot \sum_{j=0}^{r\wedge m}{d_t^{(i)} + m-1 \choose m-j}\cdot {r \choose j}\\
&= \sum_{j=0}^{m \wedge R} {d_t^{(i)} + m-1 \choose m-j} \sum_{r=j}^R
 {r \choose j}\cdot \mathds{1}\{\Delta d^{(i)}_t = r\}.\end{align*}
This gives
\begin{equation} \label{eq:sum_with_combin}
\E^{(i)}_n\left[\left. {d_t^{(i)} + \Delta d_t^{(i)}+m-1 \choose m}\right| \mathcal{F}_t\right] = \sum_{j=0}^{m \wedge R} {d_t^{(i)} + m-1 \choose m- j} \sum_{r=j}^R
 {r \choose j}\cdot \P^{(i)}_n(\Delta d^{(i)}_t = r \mid \mathcal{F}_t).
\end{equation}
We now consider the terms in the outer sum on the right-hand side  separately for different values of~$j$. For~$j = 0$, we have
$$ {d_t^{(i)} + m-1 \choose m-0} \sum_{r=0}^R
 {r \choose 0}\cdot \P^{(i)}_n(\Delta d^{(i)}_t = r \mid \mathcal{F}_t) =  {d_t^{(i)}+m-1 \choose m}.$$
For~$j = 1$,
\begin{align*}
{d_t^{(i)} + m-1 \choose m-1} \sum_{r=1}^R
 {r \choose 1}\cdot \P^{(i)}_n(\Delta d^{(i)}_t = r \mid \mathcal{F}_t) &=  {d_t^{(i)}+m-1 \choose m}\cdot \frac{m}{d_t^{(i)}} \cdot \E^{(i)}_n[\Delta d_t^{(i)}\mid \mathcal{F}_t]\\
&= {d_t^{(i)}+m-1 \choose m}\cdot \frac{m\cdot \beta_{t+1}}{t}.
\end{align*}
Where, for~$j \ge 2$, we bound:
\begin{align}
\nonumber&{d_t^{(i)} + m-1 \choose m- j} \sum_{r=j}^R
 {r \choose j}\cdot \P^{(i)}_n(\Delta d^{(i)}_t = r \mid \mathcal{F}_t)\\ 
&\le {d_t^{(i)}+m-1 \choose m}\cdot \frac{m(m-1)\cdots (m-j+1)}{d_t^{(i)}(d_t^{(i)}+1)\cdots (d_t^{(i)}+j-1)}   \cdot 2^R \cdot \P^{(i)}_n(\Delta d_t^{(i)} \ge j \mid \mathcal{F}_t).\label{eq:aux_with_2R}
\end{align}
Using the simple bound
$$\P^{(i)}_n(\Delta d_t^{(i)} \ge j \mid \mathcal{F}_t) \le 2^R \cdot \left(\frac{d_t^{(i)}}{Rt} \right)^j,$$
the right-hand side of~\eqref{eq:aux_with_2R} is smaller than
$${d_t^{(i)}+m-1 \choose m}\cdot \frac{(4m)^R}{t^j} .$$
Plugging these results in~\eqref{eq:sum_with_combin}, we obtain
\begin{align*} 
\E^{(i)}_n\left[\left. {d_t^{(i)} + \Delta d_t^{(i)}+m-1 \choose m}\right| \mathcal{F}_t\right] \le {d_t^{(i)}+m-1 \choose m}\cdot\left(1+ \frac{m\cdot \beta_{t+1}}{t} + \frac{R(4m)^R}{t^2} \right). 
\end{align*}
Hence,
$$\E^{(i)}_n[Y^{(i,m)}_{t+1} \mid \mathcal{F}_t] \le \frac{\phi_{t}^{(m)}}{\phi_{t+1}^{(m)}}\cdot Y^{(i,m)}_t \cdot \left(1+ \frac{m\cdot \beta_{t+1}}{t} + \frac{R(4m)^R}{t^2} \right) = Y^{(i,m)}_t. $$
\end{proof}

We will need the following relations between the normalization factors $\phi^{(m)}_t$ and $\phi_t$.
\begin{lemma}\label{lem:constants_k}
For every~$m \ge 2$, the limit
\begin{equation}\label{eq:phikbound}
C_m:=\lim_{t \to \infty} \frac{\phi_t^{(m)}}{(\phi_t)^m}
\end{equation}
exists, is finite and strictly positive.
\end{lemma}
\begin{proof}
For~$t\ge 2$, define
\begin{align*}
&F(t) := 1 + m\cdot \frac{\beta_{t}}{t-1},\qquad E_1(t) := \left(1+\frac{\beta_{t}}{t-1}\right)^m - F(t), \qquad E_2(t):= \frac{R(4m)^R}{(t-1)^2},
\end{align*}
so that
$$\frac{\phi(t)^k}{\phi_k(t)} = \prod_{s=1}^{t-1} \frac{F(s) + E_1(s)}{F(s)+E_2(s)} = \prod_{s=1}^{t-1}\left(1+ \frac{E_1(s) - E_2(s)}{F(s)+E_2(s)}\right),\quad t \ge 1.$$
It is easy to check that
$$F(s) + E_2(s)  \xrightarrow{s \to \infty}1,\qquad \sum_{s=1}^\infty E_1(s)< \infty,\qquad \sum_{s=1}^\infty E_2(s) < \infty.$$
Putting these observations together, the result easily follows.
\end{proof}
The following is an easy consequence of the above lemma.
\begin{lemma} \label{lem:dtikn_bounded}
For every~$m \ge 2$ there exists~$C_m' > 0$ such that
$$\mathbb{E}[(d_t^{(i)}/\phi_t)^{m}] \le C_m' \quad \text{for all } i \in \mathbb{N} \text{ and } t \in \mathbb{N}.$$
\end{lemma}
\begin{proof}
Fix~$m \ge 2$ and~$i,n\in \mathbb{N}$. For all~$t < n$ we have~$\mathbb{E}_{n}^{(i)}[(d_t^{(i)}/\phi(t))^{m}] = 0,$ 
since almost surely under~$\mathbb{P}_{n}^{(i)}$ we have~$\tau^{(i)} = n > t$, so~$d_t^{(i)} = 0$. For every~$t \ge n$ we have, for some~$C_m'$ that does not depend on~$i$,~$t$ or~$n$,
\begin{equation}
\label{eq:degconcond}\begin{split}
\E_{n}^{(i)}[       (d_t^{(i)}/\phi_t)^m       ] &\leq C_m'\cdot \E_{n}^{(i)}[        Y^{(i,m)}_t      ]  \le C_m'\cdot \mathbb{E}_{n}^{(i)}[Y_{n}^{(i,m)}] = C_m'\cdot (\phi^{(m)}_n)^{-1} \le C_m',\end{split}
\end{equation}
where the first inequality follows from Lemma~\ref{lem:constants_k} and the second inequality from Lemma~\ref{lem:Y_super}. Now,  since~$\P(\tau^{(i)}<\infty) = 1$, we have
$$\mathbb{E}[(d_t^{(i)}/\phi_t)^m] = \sum_{n=1}^\infty \mathbb{P}(\tau^{(i)}=n)\cdot \mathbb{E}_{n}^{(i)}[  (d_t^{(i)}/\phi_t)^m] \stackrel{\eqref{eq:degconcond}}{\le} C_m'.$$

\end{proof}
A consequence of the results we have so far for the process $\{X_t^{(i)}\}_{t\ge 1}$ is its convergence in the almost sure and $L_p$ senses.
\begin{proposition}
\label{prop:zetaasconv}
For each~$i \in \mathbb{N}$, the process~$\{X^{(i)}_t\}_{t \ge 1}$ converges almost surely and in~$L_q(\mathbb{P})$ (for any~$q \in [1,\infty)$) to a non-negative random variable~$\zeta^{(i)}$.
\end{proposition}
\begin{proof} By Lemma~\ref{lem:x_mart}, we have that~$(X^{(i)}_t)_{t \ge 1}$ is a non-negative submartingale. Lemma~\ref{lem:dtikn_bounded} implies that this process is bounded in~$L_m(\mathbb{P})$ for every~$m \in \mathbb{N}$; hence, it is bounded in~$L_q(\mathbb{P})$ for every~$q \in [1,\infty)$. The result then readily follows from (sub)martingale theory.
\end{proof}

Now we will establish two important properties of the random variables~$\zeta^{(i)}$ given in Proposition~\ref{prop:zetaasconv}. The first is an upper bound for its moments. This will be very useful in Section~\ref{sec:maximum}, where we study~$\max_i d_t^{(i)}$ and~$\sup_i \zeta^{(i)}$.
\begin{lemma}[Moment bounds for~$\zeta^{(i)}$]
\label{l:zetakbound}
For any~$m \ge 2$ there exists~$C''_m >0$  such that
\begin{equation}
\label{eq:zetakmoment}
\E[ (\zeta^{(i)})^m  ] \leq C''_m \cdot i^{-\beta m},\qquad i \in \mathbb{N}.
\end{equation}
\end{lemma}
\begin{proof} Fix~$m \in \mathbb{N}$. For now, also fix~$i,n \in \mathbb{N}$. We have by Proposition~\ref{prop:zetaasconv} that~$X_t^{(i)} = d_t^{(i)}/\phi(t) \xrightarrow{t \to \infty} \zeta^{(i)}$~$\mathbb{P}$-almost surely (hence also~$\mathbb{P}_{n}^{(i)}$-almost surely), so 
\[
Y_{t}^{(i,m)}=\frac{(d_{t}^{(i)}+m-1)\cdots (d_{t}^{(i)}+1)\cdot d_{t}^{(i)}}{\phi^{(m)}_t\cdot m!}\xrightarrow{t\to\infty } \frac{(\zeta^{(i)})^m}{C_m\cdot m!} \qquad\P_{n}^{(i)}\text{-a.s.},
\]
where~$C_m$ is the constant of Lemma~\ref{lem:constants_k}.
Moreover, since~$\{Y_{t}^{(i,m)}\}_{t\geq n}$ is a positive supermartingale under~$\P_{n}^{(i)}$ (by Lemma~\ref{lem:Y_super}), we have
\[
\E_{n}^{(i)}\left[\frac{(\zeta^{(i)})^m}{C_m \cdot m!} \right] \le \E_{n}^{(i)}[Y_n^{(i,m)}] = \frac{1}{\phi_n^{(m)}} \quad \Longrightarrow \quad \E_{n}^{(i)}[(\zeta^{(i)})^m] \le \frac{C_m\cdot m!}{\phi_n^{(m)}} \le \frac{C}{(\phi_n)^m}
\]
for some~$C > 0$ depending only on~$m$, by Lemma~\ref{lem:constants_k}.
 Using the fact that~$\tau(i) \ge \lfloor i/R \rfloor$, we then obtain
\[
\E[(\zeta^{(i)})^m]
=
\sum_{n=i}^{\infty} \P\left(\tau^{(i)}=n \right)\cdot \E_{n}^{(i)}[(\zeta^{(i)})^m]\leq   \frac{C}{(\phi_{\lfloor i/R \rfloor})^m}.
\]
Since~$\phi_t/t^\beta \xrightarrow{t\to \infty} b>0$, we have~$\phi_t \ge c t^\beta$ for some~$c > 0$; this completes the proof.
\end{proof}

As we discussed at the beginning of this section, another important property of the random variables~$\zeta^{(i)}$ is its positiveness. In the next result we will prove is that they are almost surely positive, which means that the process $\{d_t^{(i)}\}_{t\ge 1}$ goes to infinity with the same rate as~$\phi_t$. 
\begin{proposition}[Positiveness of~$\zeta^{(i)}$]\label{prop:zetaposi} For each~$i \in \mathbb{N}$,~$\zeta^{(i)}$ is almost surely strictly positive.
\end{proposition}
The proof of the above proposition follow ideas from \cite{S18} in the context of Balls and Bins models with immigration. The following lemma is the key ingredient needed for this result.
\begin{lemma}\label{lem:to_laplace}
Let~$i \in \mathbb{N}$. For any~$s,t \in \mathbb{N}$ with~$s < t$ we have
\begin{equation}\label{eq:main_laplace}
\mathbb{E}\left[\exp\left\{-\lambda \cdot \frac{d_t^{(i)}}{\phi_t}\right\}\right] \le \mathbb{E} \left[\exp\left\{ -\left(\lambda - \lambda^2R\cdot  \sum_{u = s}^{t-1}\frac{1}{u \cdot \phi_u} \right) \cdot \frac{d_s^{(i)}}{\phi_s} \right\} \right],\quad  \lambda > 0.
\end{equation}
\end{lemma}
We postpone the proof of this lemma to the end of this section. For now, let us see how it implies Proposition~\ref{prop:zetaposi}, which then gives Theorem~\ref{thm:convergence}.

\begin{proof}[Proof of Proposition~\ref{prop:zetaposi}]
	The proof is based on the inequality, for fixed $i$,
	\[
	\P\left( \zeta^{(i)} = 0 \right) = \P\left( \exp \{ -\lambda  \zeta^{(i)} \} = 1 \right) \le \E\left[ \exp \{ -\lambda  \zeta^{(i)} \}\right], \quad \lambda > 0,
	\]
	which follows from the Markov inequality. Our strategy is to prove, using Lemma~\ref{lem:to_laplace}, that the right-hand side can be made arbitrarily small by taking~$\lambda$ large.

Fix~$i \in \mathbb{N}$. Also fix~$s  \in \mathbb{N}$. Using Proposition~\ref{prop:zetaasconv}, dominated convergence and~\eqref{eq:main_laplace}, we have
\begin{align*}
\mathbb{E}\left[\exp\{-\lambda \zeta^{(i)}\}\right] = \lim_{t \to \infty} \mathbb{E}\left[\exp\left\{-\lambda \cdot \frac{d_t^{(i)}}{\phi_t}\right\}\right] \le \mathbb{E} \left[\exp\left\{ -\left(\lambda - \lambda^2R\cdot  \sum_{u = s}^{\infty}\frac{1}{u \cdot \phi_u} \right) \cdot \frac{d_s^{(i)}}{\phi_s} \right\} \right]
\end{align*}
for all~$\lambda > 0$. Now, using~\eqref{eq:gamma_and_b}, assume~$s$ is large enough that
$
R\sum_{u=s}^\infty \frac{1}{u\cdot \phi_u} <\frac{2}{b\beta s^{\beta}}$
and take~$\lambda = b\beta s^{\beta}/4$; then, the right-hand side above is at most~$\E\left[\exp\left \lbrace - \frac{b \beta s^{\beta}}{4}\cdot \frac{d_s^{(i)}}{\phi_s} \right \rbrace\right]$. Now, since~$\mathbb{P}(d_s^{(i)} \xrightarrow{s \to \infty} \infty) = 1$,  again using dominated convergence we obtain~$\E\left[\exp\left \lbrace - \frac{b \beta s^{\beta}}{4}\cdot \frac{d_s^{(i)}}{\phi_s} \right \rbrace\right] \xrightarrow{s\to \infty} 0$.
 This concludes the proof.
\end{proof}
Now Theorem \ref{thm:convergence} is a straightforward consequence of the results we have proven up to this point.
\begin{proof}[Proof of Theorem~\ref{thm:convergence}]
The statement readily follows from Lemma~\ref{lem:convergence_of_phi} and Propositions~\ref{prop:zetaasconv} and~\ref{prop:zetaposi}, with
$$\xi^{(i)} = \zeta^{(i)}\cdot b,$$
where~$b$ is the limit that appears in Lemma~\ref{lem:convergence_of_phi}.
\end{proof}

We finally conclude this section proving Lemma~\ref{lem:to_laplace}.
\begin{proof}[Proof of Lemma~\ref{lem:to_laplace}]
We will argue by induction on~$t-s$. Hence, for the base case of the induction, we prove that for any~$s \ge 1$,
\begin{equation}\label{eq:induction_for_laplace_base}
\mathbb{E}\left[\exp\left\{-\lambda \cdot \frac{d_{s+1}^{(i)}}{\phi_{s+1}} \right\}\right] \le \mathbb{E}\left[\exp \left\{-\left(\lambda - \lambda^2R \cdot \frac{1}{s\cdot \phi_s}\right)\cdot \frac{d_s^{(i)}}{\phi_s}  \right\} \right],\quad \lambda > 0.
\end{equation}
To prove this, note that the left-hand side is
\begin{align*}
&\mathbb{E}\left[\exp \left\{-\frac{\lambda}{\phi_{s+1}}\cdot d_s^{(i)} - \frac{\lambda}{\phi_{s+1}}\cdot \Delta d_s^{(i)}\right\}\right]\\
&\le \mathbb{E}\left[\exp\left\{- \frac{\lambda}{\phi_{s+1}}\cdot d_s^{(i)}\right\} \cdot \left(1 - \frac{\lambda}{\phi_{s+1}}\cdot \Delta d_s^{(i)} + \frac{\lambda^2}{2(\phi_{s+1})^2}\cdot (\Delta d_s^{(i)})^2 \right) \right],
\end{align*}
where the inequality  follows from applying~$e^{-x} \le 1 - x +x^2/2$, which holds for all~$x \ge 0$. Now, taking~$\mathbb{E}[\cdot \mid \mathcal{F}_s]$ inside the above expectation and using~$1-x \le e^{-x}$, the above is smaller than
\begin{equation}\label{eq:long_expon}\mathbb{E}\left[\exp\left\{- \frac{\lambda}{\phi_{s+1}}\cdot d_s^{(i)} - \frac{\lambda}{\phi_{s+1}} \cdot \mathbb{E}[\Delta d_s^{(i)}\mid \mathcal{F}_s] + \frac{\lambda^2}{2(\phi_{s+1})^2}\cdot \mathbb{E}[(\Delta d_s^{(i)})^2\mid \mathcal{F}_s]\right\}  \right].\end{equation}
Recalling that~$\mathbb{E}[\Delta d_s^{(i)}\mid \mathcal{F}_s] = \frac{\beta_{s+1}}{s}\cdot d_s^{(i)}$,
we have
\begin{equation}\begin{split}\frac{\lambda}{\phi_{s+1}}\cdot (d_s^{(i)} + \mathbb{E}[\Delta d_s^{(i)} \mid \mathcal{F}_s])& = \frac{\lambda}{\phi_{s+1}}\cdot \left(1+ \frac{\beta_{s+1}}{s}\right)\cdot d_s^{(i)} =  \frac{\lambda}{\phi_s}\cdot d_s^{(i)}. \label{eq:aux_long_expon}\end{split}\end{equation}
Next, using~$\Delta d_s(i) \in \{0,1,,\dots, R\}$, we have:
\begin{equation}\begin{split}
\mathbb{E}[(\Delta d_s^{(i)})^2\mid \mathcal{F}_s] &\le R \cdot \mathbb{E}[\Delta d_s^{(i)}\mid \mathcal{F}_s] = R\cdot \frac{\beta_{s+1}}{s}\cdot d_s^{(i)} \le \frac{R}{s}\cdot d_s^{(i)}.
\label{eq:aux_long_expon2}\end{split}\end{equation}
Now, using~\eqref{eq:aux_long_expon} and~\eqref{eq:aux_long_expon2}, the expression in~\eqref{eq:long_expon} is bounded  from above by
$$\mathbb{E}\left[\exp\left\{- \lambda \cdot \frac{d_s^{(i)}}{\phi_s} + \frac{\lambda^2 \cdot R\cdot  d_s^{(i)}}{2\cdot (\phi_{s+1})^2\cdot s}\right\}\right] \le \mathbb{E}\left[\exp\left\{-\left(\lambda - \frac{\lambda^2R}{s\cdot \phi_s}\right)\cdot \frac{d_s^{(i)}}{\phi_s} \right\}\right],$$
where the inequality follows from~$\phi_{s+1} > \phi_s$. This concludes the proof of~\eqref{eq:induction_for_laplace_base}.

We now assume that for some~$\ell \in \mathbb{N}$, we have proved the inequality in the statement of the lemma for all~$s,t$ with~$1 \le t-s \le \ell$ and all~$\lambda > 0$. Fix~$s \ge 2$ and let~$t= s+\ell$; we will carry out the induction step by proving that
\begin{equation}
\label{eq:exp_induc_step}
\mathbb{E}\left[-\lambda \cdot \frac{d_t^{(i)}}{\phi_t}\right] \le \mathbb{E}\left[\exp\left\{-\left(\lambda - \lambda^2 R \cdot \sum_{u=s-1}^{t-1} \frac{1}{u\cdot \phi_u}\right)\cdot \frac{d_{s-1}^{(i)}}{\phi_{s-1}}\right\} \right],\quad \lambda > 0.
\end{equation}
To do so, we observe that if~$\lambda > 0$ is such that the expression inside parentheses on the right-hand side is negative, then the inequality trivially holds. So we will assume from now on that
\begin{equation}\label{eq:on_lamb_poss}0 < \lambda < \left(R \cdot  \sum_{u=s-1}^{t-1}\frac{1}{u \cdot \phi(u)}\right)^{-1}.\end{equation}
By the induction hypothesis (applied to the pair~$s,t$), we have
$$\mathbb{E}\left[\exp\left\{-\lambda \cdot \frac{d_t^{(i)}}{\phi_t}\right\}\right] \le \mathbb{E}\left[ \exp\left\{- \left(\lambda - \lambda^2R   \sum_{u=s}^{t-1}\frac{1}{u\cdot \phi_u}\right)\cdot \frac{d_s^{(i)}}{\phi_s} \right\}\right].$$
Note that~\eqref{eq:on_lamb_poss} implies that the expression inside parentheses on the right-hand side is positive. So we can again apply the induction hypothesis (this time to the pair~$s-1,s$) to obtain that the right-hand side above is smaller than 
\begin{align*}&\mathbb{E}\left[\exp\left\{- \left(\lambda - \lambda^2R  \cdot \sum_{u=s}^{t-1}\frac{1}{u\cdot \phi_u} -  \left(\lambda - \lambda^2R \cdot  \sum_{u=s}^{t-1}\frac{1}{u\cdot \phi_u} \right)^2 \cdot \frac{1}{(s-1)\cdot \phi_{s-1}}\right) \cdot \frac{d_{s-1}^{(i)}}{\phi_{s-1}} \right\} \right]\\
&\le \mathbb{E}\left[\exp\left\{- \left(\lambda - \lambda^2R \cdot  \sum_{u=s}^{t-1}\frac{1}{u\cdot \phi_u} -  \lambda^2 \cdot \frac{1}{(s-1)\cdot \phi_{s-1}}\right) \cdot \frac{d_{s-1}^{(i)}}{\phi_{s-1}} \right\} \right] \\
&= \mathbb{E}\left[\exp\left\{- \left(\lambda - \lambda^2R \cdot  \sum_{u=s-1}^{t-1}\frac{1}{u\cdot \phi_u} \right) \cdot \frac{d_{s-1}^{(i)}}{\phi_{s-1}} \right\} \right]. \end{align*}
\end{proof}

\section{Central limit theorem}\label{sec:clt}
%!TEX root = ./ms.TEX
In this section, we prove Theorem~\ref{thm:clt}, which gives a central limit theorem for the process $\{X_t^{(i)} \}_{t\ge 1}$, the normalized cardinality of block $i$. For the sake of organization, we first assume Propositions \ref{prop:e_incr_d_sq} and \ref{prop:incr_d_sq} stated below and show how Theorem \ref{thm:clt} follows from them. Then, the remainder of this section is dedicated to prove both propositions.  So, throughout this section we will fix a block index $i \in \mathbb{N}$.
\subsection{Martingale central limit theorem}
Theorem \ref{thm:clt} will be obtained as a consequence of a martingale central limit theorem from~\cite{MartingaleCLTbook2} which we replicate in the appendix (Theorem~\ref{thm:clt_martingale} in Section~\ref{ss:clt_martingale}). As said above, the theorem will follow from the two propostions below.
\begin{proposition} \label{prop:e_incr_d_sq}
For any~$n$, we have that
$$\phi_t \cdot \sum_{s \ge t} \E^{(i)}_n[(\Delta X^{(i)}_s)^2] \xrightarrow{t \to \infty} \begin{cases} \E_n^{(i)}[\zeta^{(i)}]&\text{if } \beta < 1;\\[.2cm]
\E_n^{(i)} \left[\zeta^{(i)}\cdot \left( 1- \frac{\zeta^{(i)}b}{R} \right)\right]& \text{if } \beta = 1.\end{cases} $$
\end{proposition}

\begin{proposition} \label{prop:incr_d_sq}
	We have that, almost surely,
	\begin{equation*}
	\phi_t\cdot \sum_{s \ge t} (\Delta X_s^{(i)})^2 \xrightarrow{t \to \infty} \begin{cases} \zeta^{(i)}&\text{if } \beta < 1;\\ \zeta^{(i)}\cdot \left(1-\frac{\zeta^{(i)}b}{R}\right)&\text{if } \beta = 1.\end{cases}
	\end{equation*}
\end{proposition}

\begin{remark} Note that the convergence in Proposition~\ref{prop:incr_d_sq} holds~$\P$-almost surely, hence also~$\P^{(i)}_n$-almost surely for any~$n$. A natural strategy to prove these two results would be to start proving  Proposition~\ref{prop:incr_d_sq}, and then give a justification for swapping a limit with an expectation to obtain Proposition~\ref{prop:e_incr_d_sq}. This is not, however, the approach we follow: since we have not found a way to justify the swapping that is easier than computing the limit in Proposition~\ref{prop:e_incr_d_sq} directly, we just carry out the direct computation.
\end{remark}

We postpone the proofs of both propositions to Sections \ref{ss:exp_incr} and \ref{ss:sum_sq_incr} repesctively. Let us now see how they are combined in order to yield Theorem~\ref{thm:clt}.

\begin{proof}[Proof of Theorem~\ref{thm:clt}]
Fix~$n \in \N$ and define~$s_t:= \big(\sum_{s \ge t} \E_n^{(i)}[(\Delta X_s^{(i)})^2]\big)^{1/2}$ for~$t \ge n$.
Proposition~\ref{prop:e_incr_d_sq} implies that~$\sqrt{\phi_t}\cdot s_t$ converges to a positive limit, and moreover, using~$\phi_t/t^\beta \to b$ and~$\zeta^{(i)} = \xi^{(i)}/b$,
\begin{equation}\label{eq:exchange}
t^{\beta/2}\cdot s_t \xrightarrow{t \to \infty} \begin{cases}
\frac{1}{b}\cdot\left(\E_n^{(i)}[\xi^{(i)}] \right)^{1/2}&\text{if }\beta < 1;\\[.2cm]
\frac{1}{b}\cdot\left(\E_n^{(i)}\left[\xi^{(i)}\cdot \left(1-\frac{\xi^{(i)}}{R}\right)\right]\right)^{1/2}&\text{if }\beta = 1.
\end{cases}
\end{equation}
Whereas, Lemma~\ref{lem:x_mart} guarantees that~$\{X_t^{(i)}\}_{t \ge n}$ is a martingale under~$\P_n^{(i)}$;  let us show that, under this probability measure, this process satisfies the assumptions of Theorem~\ref{thm:clt_martingale}.  
We will use the bound
$$|\Delta X_t^{(i)}| \le \frac{2R}{\phi_t},\quad t \ge 1$$
obtained in~\eqref{eq:bound_triangle_z}. We have
$$
\frac{1}{s_t}\cdot \sup_{s \ge t}|\Delta X_s^{(i)}| \le \frac{1}{\sqrt{\phi_t}\cdot s_t}\cdot \sqrt{\phi_t} \cdot \frac{2R}{{\phi_t}} \xrightarrow{t \to \infty} 0
$$
and
$$\frac{1}{(s_t)^2}\cdot \E\left[\sup_{s \ge t} (\Delta X_s^{(i)})^2\right] \le \frac{1}{\phi_t\cdot (s_t)^2}\cdot \phi_t\cdot \left(\frac{2R}{\phi_t}\right)^2 \xrightarrow{t \to \infty} 0,$$
so conditions~\eqref{eq:mart_clt_2} and~\eqref{eq:mart_clt_22} are satisfied. Condition~\eqref{eq:mart_clt_222} is given by Propositions~\ref{prop:e_incr_d_sq} and~\ref{prop:incr_d_sq} together: 
$$\frac{1}{(s_t)^2} \cdot \sum_{s=t}^\infty (\Delta X_s^{(i)})^2 \xrightarrow[\P_n^{(i)}\mathrm{-a.s.}]{t \to \infty} \eta^2,$$
where
$$\eta^2 =  \frac{\xi^{(i)}}{\E_n^{(i)}[\xi^{(i)}]} \text{ if }\beta < 1\qquad \text{and}\qquad \eta^2 = \frac{\xi^{(i)}\cdot \left(1-\frac{\xi^{(i)}}{R}\right)}{\E_n^{(i)}\left[\xi^{(i)}\cdot \left(1-\frac{\xi^{(i)}}{R}\right)\right]}\text{ if }\beta = 1.$$
Hence, the conclusion of Theorem~\ref{thm:clt_martingale} tells us that, under~$\P^{(i)}_n$,
$$\frac{1}{s_t}\cdot (X_t^{(i)} - \zeta^{(i)}) = -\frac{1}{s_t} \cdot \sum_{s \ge t}\Delta X_s^{(i)} \xrightarrow[\mathrm{(d)}]{t \to \infty} \nu^{(i)}_n, $$
where~$\nu^{(i)}_n$ is the distribution of~$W \cdot \mathcal{Z}'$, where~$W,\mathcal{Z}'$ are independent,~$W$ is a standard Gaussian and the law of~$\mathcal{Z}'$ is equal to the law of~$\eta$ (under~$\P^{(i)}_n$).
Now, using~\eqref{eq:exchange},
$$t^{\beta/2}\cdot \left(\frac{d_t^{(i)}}{t^\beta} - \xi^{(i)}\right) = (t^{\beta/2}\cdot s_t \cdot b) \cdot \frac{1}{s_t} \cdot \left(\frac{\phi_t}{bt^{\beta}}\cdot \frac{d_t^{(i)}}{\phi_t} - \zeta^{(i)} \right) \xrightarrow[\mathrm{(d)}]{t \to \infty} \mu^{(i)}_n,$$
where~$\mu^{(i)}_n$ is the distribution of~$W \cdot \mathcal{Z}$, where again~$W,\mathcal{Z}$ are independent,~$W$ is a standard Gaussian and
$$\mathcal{Z}^2 \stackrel{\mathrm(d)}{=} \begin{cases}
\xi^{(i)}&\text{if } \beta<1,\\
\xi^{(i)}\cdot \left(1-\frac{\xi^{(i)}}{R}\right)&\text{if }\beta = 1,
\end{cases}$$
the distribution of the random variables on the right-hand side being under~$\P^{(i)}_n$.

The statement of the theorem now follows from the fact that~$\P(\cdot) = \sum_{n} \P^{(i)}_n(\cdot)\cdot \P(\tau^{(i)} = n)$, so the law~$\mu^{(i)}$ defined there is equal to~$\sum_n \mu^{(i)}_n\cdot \P(\tau^{(i)}=n)$.
\end{proof}

\subsection{Sum of expected squared increments: Proof of Proposition~\ref{prop:e_incr_d_sq}}\label{ss:exp_incr}
We now give the proof of our statement concerning the asymptotic behavior of the sum~$\sum_{s \ge t}\mathbb{E}[(\Delta X^{(i)}_s)^2]$. For this proof and several of the following, the equality~$\Delta X^{(i)}_t = A_t - B_t$ from~\eqref{eq:a_and_b} will play an important role, for this reason we recall $A_t$ and $B_t$ below 
\begin{equation*}
    A_t:= \frac{\Delta d^{(i)}_t}{\phi_{t+1}},\qquad B_t:= X^{(i)}_t \cdot \frac{\Delta \phi_t}{\phi_{t+1}} = X^{(i)}_t \cdot \left(\frac{\phi_{t+1}}{\phi_t} - 1\right) = X^{(i)}_t \cdot \frac{\beta_{t+1}}{t}.
    \end{equation*}
We also recall, for future reference, that
\begin{equation} \label{eq:asymp_of_phi}
\frac{\phi_t}{t^\beta}\xrightarrow{t \to \infty} b,\quad \frac{\phi_t}{\phi_{t+1}} \xrightarrow{t \to \infty} 1,\quad \frac{\Delta \phi_t}{\phi_t} = \frac{\phi_{t+1}}{\phi_t}- 1 = \frac{\beta_t}{t},\quad \quad \beta_t \xrightarrow{t \to \infty} \beta.
\end{equation}

\begin{proof}[Proof of Proposition~\ref{prop:e_incr_d_sq}]
As in Lemma~\ref{lem:expectation_delta_d}, we have
$$\E^{(i)}_n[\Delta d^{(i)}_t \mid \mathcal{F}_t] = \sum_{r=0}^R g_{t+1}(r)\cdot r \cdot \frac{d_t^{(i)}}{Rt} = \beta_{t+1} \cdot \frac{d^{(i)}_t}{t},\qquad t \ge n. $$
By a similar proof as the one of that lemma, using the formula for the second moment of binomial random variables, we also get, for~$t \ge n$,
$$\E^{(i)}_n[(\Delta d^{(i)}_t)^2 \mid \mathcal{F}_t] = \sum_{r=0}^R g_{t+1}(r) \cdot \left(r\cdot \frac{d_t^{(i)}}{Rt} + (r^2-r)\cdot \left(\frac{d_t^{(i)}}{Rt} \right)^2 \right) = \beta_{t+1}\cdot \frac{d_t^{(i)}}{t} + \beta_{t+1}' \cdot \left(\frac{d_t^{(i)}}{t} \right)^2,$$
\iffalse
We start noting that
$$\P(\Delta d^{(i)}_t = r\mid \mathcal{F}_t) = \sum_{u=0}^R g_t(u)\cdot \mathscr{B}\left(u,\frac{d^{(i)}_t}{Rt},r\right),\qquad r \in \{0,\ldots, R\},$$
that is, the law of~$\Delta d^{(i)}_t$ conditioned on~$\mathcal{F}_t$ is a mixture of binomials: for~$0 \le u \le R$, the binomial distribution with parameters~$u$ and~$\frac{d^{(i)}_t}{Rt}$ has weight~$g_t(u)$. The formulas for first and second moments for the binomial distribution then give:
\begin{align*}
&\E[\Delta d^{(i)}_t \mid \mathcal{F}_t] = \sum_{r=0}^R g_t(r)\cdot r \cdot \frac{d_t^{(i)}}{Rt} = \beta_t \cdot \frac{d^{(i)}_t}{t},\\[.2cm]
&\E[(\Delta d^{(i)}_t)^2 \mid \mathcal{F}_t] = \sum_{r=0}^R g_t(r) \cdot \left(r\cdot \frac{d_t^{(i)}}{Rt} + (r^2-r)\cdot \left(\frac{d_t^{(i)}}{Rt} \right)^2 \right) = \beta_t \cdot \frac{d_t^{(i)}}{t} + \beta_t' \cdot \left(\frac{d_t^{(i)}}{t} \right)^2,
\end{align*}
\fi
where we define~$\beta_t':= R^{-2}\sum_{r=0}^R g_t(r)\cdot (r^2-r)$. Using this together with~\eqref{eq:a_and_b} we obtain
\begin{align*}
\E^{(i)}_n[(\Delta X^{(i)}_t)^2\mid \mathcal{F}_t] &= \E^{(i)}_n[(A_t)^2 + (B_t)^2 - 2A_tB_t \mid \mathcal{F}_t]\\
&=  \frac{\beta_{t+1}\cdot d_t^{(i)}}{(\phi_{t+1})^2\cdot t} + \frac{\beta'_{t+1} \cdot (d_t^{(i)})^2}{(\phi_{t+1})^2\cdot t^2} + \frac{(X_t^{(i)})^2 \cdot (\Delta \phi_t)^2}{(\phi_{t+1})^2} - 2 \frac{\beta_{t+1} \cdot X_t^{(i)}\cdot \Delta \phi_t \cdot d_t^{(i)}}{(\phi_{t+1})^2 \cdot t},
\end{align*}
so, taking the expectation and using~$X^{(i)}_t = d^{(i)}_t/\phi_t$,
\begin{equation}\begin{split}
\E^{(i)}_n[(\Delta X^{(i)}_t)^2] = & \frac{\beta_{t+1} \cdot \phi_t}{(\phi_{t+1})^2 \cdot t}\cdot \E^{(i)}_n[X^{(i)}_t] + \left(\frac{\beta'_{t+1} \cdot (\phi_t)^2}{(\phi_{t+1})^2  \cdot t^2} + \left( \frac{\Delta \phi_t}{\phi_{t+1}}\right)^2 - 2\frac{\beta_{t+1} \cdot \Delta \phi_t \cdot \phi_t}{(\phi_{t+1})^2 \cdot t} \right)  \cdot \E^{(i)}_n[(X^{(i)}_t)^2].\end{split} \label{eq:split_terms}
\end{equation}
We let
$$\beta':= \frac{1}{R^2}\sum_{r=0}^R g_\infty(r)\cdot (r^2-r) = \lim_{t \to \infty} \beta'_t.$$
We then obtain, using~\eqref{eq:asymp_of_phi}, the following asymptotic expressions for the quotients that appear in~\eqref{eq:split_terms}:
\begin{align*}
\frac{\beta_{t+1} \cdot \phi_t}{(\phi_{t+1})^2 \cdot t} \sim \frac{\beta}{bt^{1+\beta}},\quad \frac{\beta'_{t+1} \cdot (\phi_t)^2}{(\phi_{t+1})^2  \cdot t^2} \sim \frac{\beta'}{t^2}, \quad \left( \frac{\Delta \phi_t}{\phi_{t+1}}\right)^2 \text{ and }\; \frac{\beta_{t+1} \cdot \Delta \phi_t \cdot \phi_t}{(\phi_{t+1})^2 \cdot t} \sim \left(\frac{\beta}{t}\right)^2.
\end{align*}
Recall that~$X^{(i)}_t \to \zeta^{(i)}$ in~$L_p(\P)$ for all~$p \in [1,\infty)$, hence also in~$L_p(\P^{(i)}_n)$ for all~$p \in [1,\infty)$, so
$$\E^{(i)}_n[X_t^{(i)}]\xrightarrow{t \to \infty} \E^{(i)}_n[\zeta^{(i)}]\quad \text{and} \quad \E^{(i)}_n[(X^{(i)}_t)^2] \xrightarrow{t \to  \infty } \E^{(i)}_n[(\zeta^{(i)})^2].$$

We study the behavior of~$\phi_t \cdot \sum_{s \ge t}\E^{(i)}_n[(\Delta X_s)^2]$ when~$t \to \infty$ by separately considering terms obtained from~\eqref{eq:split_terms}, using the asymptotic expressions obtained above. First,
\begin{equation}\phi_t \cdot \sum_{s \ge t} \frac{\beta_{s+1} \cdot \phi_s}{(\phi_{s+1})^2 \cdot t}\cdot \E^{(i)}_n[X^{(i)}_s] \sim  bt^\beta \cdot \sum_{s \ge t}\frac{\beta}{bs^{1+\beta}}\cdot \E^{(i)}_n[\zeta^{(i)}] \xrightarrow{t \to \infty}  \E^{(i)}_n[\zeta^{(i)}]. \label{eq:1cl}\end{equation}
Next, we have:
\begin{align} 
\label{eq:2cl}&\phi_t \cdot \sum_{s \ge t} \frac{\beta'_{s+1} \cdot (\phi_s)^2}{(\phi_{s+1})^2 \cdot s^2}\cdot \E^{(i)}_n[(X^{(i)}_t)^2] \sim \beta'b\cdot \E^{(i)}_n[(\zeta^{(i)})^2]\cdot t^{\beta-1},\\[.2cm]
\label{eq:3cl}&\phi_t \cdot \sum_{s \ge t} \left(\frac{\Delta \phi_t}{\phi_{t+1}}\right)^2 \cdot \E^{(i)}_n[(X^{(i)}_t)^2]\sim \beta^2 b \cdot \E^{(i)}_n[(\zeta^{(i)})^2]\cdot t^{\beta - 1},\\[.2cm]
\label{eq:4cl}&\phi_t \cdot \sum_{s \ge t} \frac{\beta_{t+1} \cdot \Delta \phi_t \cdot \phi_t}{(\phi_{t+1})^2 \cdot t} \cdot \E^{(i)}_n[(X^{(i)}_t)^2]\sim \beta^2 b \cdot \E^{(i)}_n[(\zeta^{(i)})^2]\cdot t^{\beta - 1}.
\end{align}
Now, when~$\beta < 1$, the desired result already follows from putting together~\eqref{eq:1cl}-\eqref{eq:4cl} in~\eqref{eq:split_terms}. For the case~$\beta =1$, the result follows in the same way, with the additional observation that in this case we have~$g_\infty(R) = 1$, so~$\beta' = 1-\frac{1}{R}$.
\end{proof}

\subsection{Sum of squared increments: Proof of Proposition~\ref{prop:incr_d_sq}} \label{ss:sum_sq_incr}
Recall that Proposition \ref{prop:incr_d_sq} concerns almost sure convergence of the random series $
\phi_t\cdot \sum_{s \ge t} (\Delta X_s^{(i)})^2$. Notice that the order of magnitude of such random series is intrinsically connected to those steps when block $i$ has increased its cardinality as well as by how much said cardinality increased. This forces us to keep track on the times block $i$ received exactly $r$ elements and investigate the behavior of $\phi_t\cdot \sum_{s \ge t} (\Delta X_s^{(i)})^2$ along these random sets. For this reason, the need for results on sums along subsets of $\mathbb{N}$ arises naturally in this section.

For the sake of organization we will prove Proposition~\ref{prop:incr_d_sq} separately for~$\beta = 1$ and~$\beta < 1$; the latter will be somewhat more involved. In order to set up the stage for the first case~$\beta = 1$, we start with the following.
\begin{definition} Define the random sets
\begin{equation*}
\mathscr{T}^{(i)}(r):= \{t: \Delta d^{(i)}_t = r \},\quad r \in \{0,\ldots, R\}.
\end{equation*}
\end{definition}
For the following statement, recall that a set~$\Lambda \subset \mathbb{N}$ has \emph{asymptotic density}~$\alpha \in [0,1]$ if
$$\lim_{N \to \infty} \frac{|\Lambda \cap \{1,\ldots, N\}|}{N} = \alpha.$$
Our first result concerns the asymptotic density of $\mathscr{T}^{(i)}(r)$ in the case $\beta =1$.
\begin{lemma}\label{lem:asymp_dens_g_infty_1}
	Assume~$\beta =1$. Then, almost surely for each~$r \in \{0,\ldots, R\}$, the set~$\mathscr{T}^{(i)}(r)$ has asymptotic  density~$\mathscr{B}\left(R, \zeta^{(i)} b/R, r\right).$
\end{lemma}
\begin{proof}
	When~$\beta = 1$, we have~$g_\infty(R) = 1$ and~$\frac{d_t^{(i)}}{t} \xrightarrow{t \to \infty} \zeta^{(i)}\cdot b$ almost surely, so for each~$r$,
	$$\mathbb{P}(\Delta d_t^{(i)} = r\mid \mathscr{F}_t) = \sum_{u=0}^R g_{t+1}(u)\cdot \mathscr{B}(u,d_t^{(i)}/(Rt),r) \xrightarrow{t \to \infty}  \mathscr{B}(R,\zeta^{(i)}b/R,r)$$
	almost surely, and then
	$$\frac{1}{t} \sum_{s=1}^t    \mathbb{P}(\Delta d_s^{(i)} = r\mid \mathcal{F}_s) \xrightarrow{t \to \infty}  \mathscr{B}(R,\zeta^{(i)}b/R,r)$$
	almost surely. The result now follows from applying the Azuma-Hoeffding inequality (Theorem~\ref{thm:azuma}) and the Borel-Cantelli lemma to the martingale
	$$\sum_{s=\new{1}}^{t} \left( \mathds{1}\{\Delta d_s^{(i)} = r\} -  \mathbb{P}(\Delta d_s^{(i)} = r\mid \mathcal{F}_s)\right),\quad t \in \mathbb{N}.$$
\end{proof}

We now state a result about the asymptotic value of series that are taken over sets with a given asymptotic density. In order to keep the flow of the presentation, we give the proof in an appendix.
\begin{lemma} \label{lem:sum_over}
	Let~$\Lambda \subset  \mathbb{N}$ be a set with asymptotic density equal to~$\alpha \in [0,1]$. Then,
	$$\lim_{N \to \infty} N \sum_{\substack{n \ge N,\\n \in \Lambda}} \frac{1}{n^2} = \alpha.$$
\end{lemma}
Putting the two previous results together, we obtain the following.
\begin{corollary} \label{lem:new_cor_lemmas}
	Assume~$\beta = 1$. Then, almost surely,
	\begin{equation}\label{eq:new_cor_1}\lim_{t \to \infty} t  \sum_{s \ge t} \frac{\Delta d_s^{(i)}}{s^2} = \zeta^{(i)}b \end{equation}
	and
	\begin{equation} \label{eq:new_cor_2} \lim_{t \to \infty} t  \sum_{s \ge t} \frac{(\Delta d_s^{(i)})^2}{s^2} = (\zeta^{(i)})^2b^2 + \zeta^{(i)}b -\frac{(\zeta^{(i)})^2b^2}{R}.\end{equation}
\end{corollary}
\begin{proof}
	For any~$h:\mathbb{R} \to \mathbb{R}$ we have
	\begin{align*}
	\sum_{s \ge t}\frac{h(\Delta d^{(i)}_s)}{s^2} = \sum_{r=0}^R h(r) \sum_{s \ge t} \frac{\mathds{1}{\{s\in \mathscr{T}^{(i)}(r)\}}}{s^2},
	\end{align*}
	so, by Lemma~\ref{lem:asymp_dens_g_infty_1} and Lemma~\ref{lem:sum_over}, we get
	\begin{equation}
	\label{eq:third_with_asymp}
	\lim_{t \to \infty} t \sum_{s \ge t} \frac{h(\Delta d^{(i)}_s)}{s^2} = \sum_{r = 0}^R h(r)\cdot \mathscr{B}(R,\zeta^{(i)}b/R,r).
	\end{equation}
	The result now follows by using the formulas for the first and second moments  of the binomial distribution.
\end{proof}

We are now ready  to treat the case~$\beta = 1$ in Proposition~\ref{prop:incr_d_sq}.
\begin{proof}[Proof of Proposition~\ref{prop:incr_d_sq}, case~$\beta = 1$] We write~$\Delta X^{(i)}_t =  A_t - B_t $ as in~\eqref{eq:a_and_b}.
\iffalse	Note that
	$$X_t^{(i)} = \frac{d_t^{(i)}}{\phi_t} \xrightarrow{t \to \infty} \zeta^{(i)},\quad \frac{\phi_t}{t} \xrightarrow{t\to \infty} b,\quad \frac{d_t^{(i)}}{t} \xrightarrow{t \to \infty} \zeta^{(i)}b,$$
	and (using~$\phi_t \xrightarrow{t\to \infty} \infty$ and~$\phi_{t+1}-\phi_t \le 1$),
	\begin{equation} \label{eq:with_phi}\frac{\phi_{t+1}}{\phi_t}\xrightarrow{t\to\infty}1 \quad \text{and}\quad\frac{\phi_{t+1}-\phi_t}{\phi_{t+1}} = \frac{\phi_t}{\phi_{t+1}}\cdot \left(\frac{\phi_{t+1}}{\phi_t} - 1\right) \sim \frac{\beta}{t} = \frac{1}{t}.\end{equation}
	By these considerations, we get
	\begin{equation*}
	A_t \sim \frac{\Delta d_t^{(i)}}{bt},\qquad  B_t \sim \frac{\zeta^{(i)}}{t}.
	\end{equation*}
\fi
We then have
	$$\phi_t \cdot \sum_{s \ge t} (\Delta X_s^{(i)})^2 = \phi_t \cdot \sum_{s \ge t} (A_s)^2 + \phi_t \cdot \sum_{s \ge t} ( B_s)^2 - 2\phi_t\sum_{s \ge t} ( A_s \cdot B_s)$$
	and we compute separately the limit of the three terms on the right-hand side as~$t \to \infty$. For the first term, we have
	$$\lim_{t \to \infty} \phi_t \sum_{s \ge t} (A_s)^2 = \lim_{t \to \infty} bt \sum_{s \ge t} \left(\frac{\Delta d_s^{(i)}}{bs} \right)^2 = \zeta^{(i)} - \frac{(\zeta^{(i)})^2b}{R} + (\zeta^{(i)})^2b, $$
	where the second equality follows from~\eqref{eq:new_cor_2}. For the second term,
	$$\lim_{t \to \infty} \phi_t \sum_{s \ge t} (B_s)^2 = \lim_{t \to \infty} (\zeta^{(i)})^2bt \sum_{s \ge t}\frac{1}{s^2} = (\zeta^{(i)})^2b. $$
	Finally, for the third term,
	$$\lim_{t \to \infty} \phi_t \sum_{s \ge t} (A_s \cdot B_s) = \lim_{t \to \infty} \zeta^{(i)}t \sum_{s \ge t}\frac{\Delta d_s^{(i)}}{s^2} = (\zeta^{(i)})^2b,$$
	where the second equality follows from~\eqref{eq:new_cor_1}. Putting things together, we then obtain
	$$\lim_{t \to \infty} \phi_t \sum_{s \ge t} (\Delta X_s)^2 = \zeta^{(i)} \left( 1- \frac{\zeta^{(i)}b}{R}\right).$$
\end{proof}

We now  give some additional definitions and preliminary results concerning both the cases~$\beta  = 1$ and~$\beta < 1$. For the missing part of the proof of Proposition~\ref{prop:incr_d_sq}, we will need the case~$\beta < 1$, but in the next section we will also use the results obtained here for~$\beta = 1$.

\begin{definition}
We define a sequence of stopping times~$(\sigma^{(i)}_n)_{n \in \mathbb{N}_0}$ by letting~$\sigma_0^{(i)}:= \tau^{(i)}$ and, for~$k \in \mathbb{N}_0$,
$$\sigma_{k+1}^{(i)}:= \inf\left\{t > \sigma_k^{(i)}: d^{(i)}_t > d^{(i)}_{\sigma_k^{(i)}}\right\} = 1+\inf\{t \ge \sigma_k^{(i)}: \Delta d_t^{(i)} \neq 0\}.$$ 
In words, $\sigma_k^{(i)}$ is the $k$-th time the process has increased the cardinality of block $i$.
We also define
$$D_k^{(i)} := d^{(i)}_{\sigma_k^{(i)}},\quad\mathcal{G}_k^{(i)} = \mathcal{F}_{\sigma_k^{(i)}},\quad k \in \mathbb{N}_0.$$
Finally, writing~$\Delta D_k^{(i)} := D_{k+1}^{(i)} - D_k^{(i)}$, define the random sets
$$\mathscr{K}^{(i)}(r):= \{k \in \mathbb{N}_0: \Delta D^{(i)}_k = r\},\qquad r \in \{1,\ldots, R\},$$
\end{definition}
That is, $\mathscr{K}^{(i)}(r)$ is the set of moments when the cardinality of block $i$ increased by exactly $r$. Note that the definition of~$\mathscr{K}^{(i)}(r)$ above does not include~$r = 0$ because by definition,~$\Delta D^{(i)}_k$ is never zero. 

The next lemma relates the asymptotic density of $\mathscr{K}^{(i)}(r)$ and $\zeta^{(i)}$.
\begin{lemma}\label{lem:asymp_dens_K}
Almost surely for each~$r \in \{1,\ldots, R\}$, the set~$\mathscr{K}^{(i)}(r)$ has asymptotic density $$\frac{\mathscr{B}(R,\zeta^{(i)}b/R,r)}{1-\mathscr{B}(R,\zeta^{(i)}b/R,0)}$$ if~$\beta = 1$ and~$\mathds{1}\{r = 1\}$ if~$\beta < 1$.
\end{lemma}
We postpone the proof of this lemma. Let us now give a second result concerning sums of series with asymptotic densities.
\begin{lemma}\label{lem:sop_ad}
	Let~$(a_k)_{k \in \mathbb{N}}$ be an increasing sequence of natural numbers with~$\Delta a_k \in \{1,\ldots, R\}$ for each~$k$, and such that for each~$r \in \{1,\ldots, R\}$, the set~$\{k:\Delta a_k = r\}$ has asymptotic density~$\rho_r$. We then have
	\begin{equation*}
	\lim_{k_0\to \infty} a_{k_0}\cdot \sum_{k \ge k_0} \left(\frac{\Delta a_k}{a_k}\right)^2 = \frac{\sum_r r^2\rho_r}{\sum_r r \rho_r }.
	\end{equation*}
\end{lemma}
Again, the proof is carried out in the Appendix. We will now obtain the following by combining the two previous lemmas:
\begin{corollary} \label{cor:g_inf_less}
	Assume~$\beta < 1$. Then,
	\begin{equation}\label{eq:of_cor2}
	\lim_{t \to \infty} \phi_t \cdot \sum_{s \ge t} \left(\frac{\Delta d^{(i)}_s}{d^{(i)}_s} \right)^2 = \frac{1}{\zeta^{(i)}}.
	\end{equation}
\end{corollary}
\begin{proof}
	Define
	$$K(t):= \sup\{k: \sigma_k^{(i)} \le t\},\quad t \in \mathbb{N};$$
	note that~$d_t^{(i)} = d_{\sigma_{K(t)}^{(i)}}^{(i)} = D^{(i)}_{K(t)}$ for each~$t$. Moreover, it is easy to check that
	$$\sum_{s \ge t} \left(\frac{\Delta d^{(i)}_s}{d^{(i)}_s} \right)^2 =  \sum_{k \ge K(t)} \left(\frac{\Delta D_k^{(i)}}{D_k^{(i)}} \right)^2 $$
	for each~$t$. Using the fact, given in Lemma~\ref{lem:asymp_dens_K}, that~$\mathscr{K}^{(i)}(r)$ has asymptotic density equal to $\mathds{1}\{r = 1\}$, together with Lemma~\ref{lem:sop_ad}, we obtain
	$$D^{(i)}_{K(t)} \sum_{k \ge K(t)} \left(\frac{\Delta D_k^{(i)}}{D_k^{(i)}} \right)^2 \xrightarrow{t \to \infty} 1, $$
	so
	$$\lim_{t \to \infty} \phi_t \cdot \sum_{s \ge t}  \left(\frac{\Delta d^{(i)}_s}{d^{(i)}_s} \right)^2 = \lim_{t \to \infty} \frac{\phi_t}{d_t^{(i)}}\cdot \lim_{t \to \infty} \left( D^{(i)}_{K(t)}  \sum_{k \ge K(t)} \left( \frac{\Delta D_k^{(i)}}{D^{(i)}_k } \right)^2 \right) = \frac{1}{\zeta^{(i)}}.$$
\end{proof}
Now we are finally able to cover the case $\beta < 1$.
\begin{proof}[Proof of Proposition~\ref{prop:incr_d_sq}, case~$\beta < 1$] As in the proof of the other case, we write $\Delta X^{(i)}_t = A_t - B_t$, with~$A_t$ and~$B_t$ defined in~\eqref{eq:a_and_b}. Again we write
	$$\phi_t \cdot \sum_{s \ge t} (\Delta X_s^{(i)})^2 = \phi_t \cdot \sum_{s \ge t} (A_s)^2 + \phi_t \cdot \sum_{s \ge t} ( B_s)^2 - 2\phi_t\sum_{s \ge t} ( A_s \cdot B_s)$$
	and consider the three terms on the right-hand side separately. For the first term,  Corollary~\ref{cor:g_inf_less} gives
	$$\lim_{t \to \infty} \phi_t \sum_{s \ge t} (A_s)^2 = \zeta^{(i)}.$$
	Next, since~$B_t \sim \zeta^{(i)} \cdot \frac{\beta}{t}$ and~$\phi_t = o(t)$ when~$\beta < 1$, we obtain
	$$\lim_{t \to \infty} \phi_t \sum_{s \ge t} (B_s)^2 = 0.$$
	Finally, by the Cauchy-Schwarz inequality,
	$$\phi_t \sum_{s \ge t} (A_s \cdot B_s) \le \sqrt{\phi_t \sum_{s \ge t}(A_s)^2 } \cdot \sqrt{\phi_t \sum_{s \ge t}(B_s)^2 } \xrightarrow{t \to \infty} 0.$$
	Putting things together, we obtain~$\phi_t \cdot \sum_{s \ge t} (\Delta X_s^{(i)})^2 \xrightarrow{t \to \infty} \zeta^{(i)}$.
\end{proof}

It remains to prove Lemma~\ref{lem:asymp_dens_K}. The proof will involve further definitions and lemmas, some of which will be useful in the next section.
\begin{definition}\label{def:P_and_P_hat} Define, for~$m, t \in \mathbb{N}$,
\begin{equation*}\begin{split}
&\hat{P}_{m,t}(r) :=  \sum_{u=0}^R g_{t+1}(u)\cdot \mathscr{B}\left(u,\frac{m}{Rt},r\right), \quad r \in \{0,\ldots, R\}. 
\end{split}
\end{equation*}
Notice that in the event $\{d_t^{(i)} = m\}$, the following identity holds
$
\hat{P}_{m,t}(r)=\mathbb{P}(\Delta d^{(i)}_t = r \mid \mathcal{F}_t).
$

We also define
$$\hat{\mathscr{P}}_{m,t}(r) := \hat{P}_{m,t}(r) + \sum_{t'=t+1}^\infty \left(\prod_{s=t}^{t'-1} \hat{P}_{m,s}(0) \right)\cdot \hat{P}_{m,t'}(r), \quad r \in \{1,\ldots, R\}.$$
Finally, also define, for~$\alpha \in (0,1]$,
$$\mathscr{P}_\alpha(r):= \frac{\mathscr{B}(R,\alpha,r)}{1-\mathscr{B}(R,\alpha,0)}, \quad r \in \{1,\ldots, R\}$$
and finally, let
$$\mathscr{P}_0(r):= \mathds{1}\{r = 1\},\quad r \in \{1,\ldots, R\}.$$
\end{definition}
Note that~$\hat{\mathscr{P}}_{m,t}(r)$ is the probability that, given a block of cardinality~$m$ at time~$t$, it receives exactly~$r$ elements at the next time it receives new elements. Since no block ever permanently stops receiving new numbers, we have that
$$\sum_{r = 1}^R \hat{\mathscr{P}}_{m,t}(r) = 1.$$
Moreover, we have
\begin{equation}\mathbb{P}(\Delta D^{(i)}_k = r \mid \mathcal{G}_k^{(i)}) = \hat{\mathscr{P}}_{D^{(i)}_k,\sigma^{(i)}_k}(r), \qquad r \in \{1,\ldots, R\}.\label{eq:key_eq_hat_P}\end{equation}
 With the above notation in mind, the next lemmas will play important roles in the proof of Lemma~\ref{lem:asymp_dens_K}.
\begin{lemma}\label{lem:conv_p_1}
Assume that~$\beta = 1$.	Let~$(m_n)$,~$(t_n)$ be increasing sequences of natural numbers with $m_n/t_n \xrightarrow{n \to \infty} \alpha R$, with~$\alpha > 0$. Then,
	$$\hat{\mathscr{P}}_{m_n,t_n}(r) \xrightarrow{n \to \infty} {\mathscr{P}}_\alpha(r) \quad \text{for all } r \in \{1,\ldots, R\}.$$
\end{lemma}
\begin{proof}
	For each~$n \in \mathbb{N}$, let
	$$\pi_n(0,r):= \hat{P}_{m_n,t_n}(r),\qquad r \in \{1,\ldots,R\}$$
	and
	$$\pi_n(t,r):=\left(\prod_{s=0}^{t-1} \hat{P}_{m_n,t_n+s}(0) \right)\cdot \hat{P}_{m_n,t_n+t}(r),\qquad t \ge 1,\;r \in \{1,\ldots, R\}.$$
	Next, define
	$$P_\alpha(r) := \mathscr{B}(R,\alpha,r),\qquad \alpha \in (0,1],\; r \in \{0,\ldots, R\}$$
	and
	$$\pi(t,r) := ({P}_\alpha(0))^{t}\cdot {P}_\alpha(r),\qquad t \ge 0,\;r \in \{1,\ldots, R\}.$$
	Now,~$\pi_n$ and~$\pi$ are probabilities on~$\mathbb{N}_0 \times\{1,\ldots R\}$ and, using the assumptions that~$\beta = 1$ (so~$g_\infty(R) = 1$) and~$m_n/t_n \to \alpha/R$ together with the definitions of~$\hat{P}_{m_n,t_n}$ and~${P}_\alpha$, it is readily seen that~$\pi_n(t,r) \to \pi(t,r)$ for every~$(t,r)$. It then follows from Scheff\'e's theorem that, for all~$r \in \{1,\ldots, R\}$, we have
	$$\hat{\mathscr{P}}_{m_n,t_n}(r) = \sum_{t =0}^\infty \pi_n(t,r) \xrightarrow{n \to \infty} \sum_{t =0}^\infty \pi(t,r) = {\mathscr{P}}_\alpha(r).$$
\end{proof}
\begin{lemma} \label{lem:conv_p_0}
	Let~$(m_n)$,~$(t_n)$ be increasing sequences of natural numbers with $m_n/t_n \xrightarrow{n \to \infty} 0$. Then,
	$$\hat{\mathscr{P}}_{m_n,t_n}(r) \xrightarrow{n \to \infty} {\mathscr{P}}_0(r) \quad \text{for all } r \in \{1,\ldots, R\}.$$
\end{lemma}
\begin{proof}
	We first claim that for any~$\varepsilon > 0$ there exists~$\delta_\varepsilon > 0$ such that
	\begin{equation*}
	0\le p < \delta_\varepsilon,\; r \in \{0,\ldots, R\},\;u\in \{2,\ldots, R\}\quad \Longrightarrow \quad \mathscr{B}(r,p,u) \le \frac{\varepsilon}{R}\cdot \mathscr{B}(r,p,1).
	\end{equation*}
	This is easily checked by first treating the  cases~$r = 0$ and~$r=1$ first (for both these cases the statement is trivial), and, for~$r \ge 2$, using the facts that~$\mathscr{B}(r,p,1) \sim rp$ as~$p \to 0$ for~$r \ge 1$, and~$\mathscr{B}(r,p,u) = o(p)$ as~$p \to 0$ for~$u \ge 2$.
	
	Now, fix~$\varepsilon > 0$ and fix~$n_0$ such that~$\frac{m_n}{Rt_n} < \delta_\varepsilon$ for all~$n \ge n_0$. Then, for all~$n \ge n_0$,~$u \in \{2,\ldots, R\}$ and~$t \ge t_n$ we have
	\begin{align*}
	\hat{P}_{m_n,t}(u) &= \sum_{r=0}^R g_{t+1}(r)\cdot \mathscr{B}\left(r,\frac{m_n}{Rt_n},u\right)\le \frac{\varepsilon}{R}\sum_{r=0}^R g_{t+1}(r)\cdot \mathscr{B}\left( r,\frac{m_n}{Rt},1\right) = \frac{\varepsilon}{R}\cdot \hat{P}_{m_n,t}(1).
	\end{align*}
	This readily gives, for all~$n \ge n_0$ and~$t \ge t_n$,
	$$\sum_{u=2}^R \hat{P}_{m_n,t}(u) \le \varepsilon\cdot  \hat{P}_{m_n,t}(1).$$
	Using the definition of~$\hat{\mathscr{P}}_{m_n,t_n}$, we then obtain
	$$\sum_{u=2}^R \hat{\mathscr{P}}_{m_n,t_n}(u) \le \varepsilon \cdot \hat{\mathscr{P}}_{m_n,t_n}(1).$$
	Combined with~$\sum_{u=1}^R \hat{\mathscr{P}}_{m_n,t_n}(u) = 1$, this gives~$1 \ge \hat{\mathscr{P}}_{m_n,t_n}(1) \ge \frac{1}{1+\varepsilon},$
	completing the proof.
\end{proof}

\begin{corollary}\label{cor:conv_of_scrP}
We have that, almost surely, for each~$r \in \{1,\ldots, R\}$,
\begin{equation*}
\P(\Delta D^{(i)}_k = r \mid \mathcal{G}_k^{(i)}) \xrightarrow{k \to \infty} \begin{cases}
\mathscr{P}_0(r)&\text{if }\beta < 1;\\[.2cm] \mathscr{P}_{\zeta^{(i)}b/R}(r)&\text{if } \beta = 1.
\end{cases}
\end{equation*}
\end{corollary}
\begin{proof}
Recall~\eqref{eq:key_eq_hat_P} and note that
	$$\frac{d^{(i)}_t}{t} \xrightarrow{t \to \infty} \begin{cases} 0&\text{if } \beta < 1;\\ \zeta^{(i)}b& \text{if }\beta = 1,\end{cases}\qquad \text{so} \qquad \frac{D^{(i)}_k}{\sigma^{(i)}_k} = \frac{d^{(i)}_{\sigma^{(i)}_k}}{\sigma^{(i)}_k} \xrightarrow{k \to \infty} \begin{cases} 0&\text{if } \beta < 1;\\ \zeta^{(i)}b& \text{if }\beta = 1.\end{cases} $$
	Hence, the desired convergence follows from Lemma~\ref{lem:conv_p_1} and Lemma~\ref{lem:conv_p_0}.
\end{proof}
\begin{proof}[Proof of Lemma~\ref{lem:asymp_dens_K}]
The desired result follows from Corollary~\ref{cor:conv_of_scrP} and a simple application of the Azuma-Hoeffding inequality to the martingale
	$$\sum_{\ell=1}^k \left(\mathds{1}\{\Delta D^{(i)}_\ell = r\}  - \P(\Delta D^{(i)}_\ell = r \mid \mathcal{G}_\ell^{(i)})\right),\quad k \in \N;$$
we omit the details.
\end{proof}

\section{Persistent leadership}\label{sec:leader}
%!TEX root = ./ms.tex
 We note that, if we had~$\P(\xi^{(i)} = \xi^{(j)}) = 0$ whenever~$i \neq j$, we would already have persistent leadership. Indeed, as we will see in the proof of Theorem~\ref{thm:leader}, the exponential decay of the tails of the normalized cardinality given by Lemma~\ref{l:upperbounddeg} implies that, almost surely, only a finite number of blocks can ``compete'' for the leadership. If, $\xi^{(i)} \neq \xi^{(j)}$ for $i \neq j$, then the cardinalities of the respective blocks at large time~$t$ must be at distance of order~$t^\beta$ from one another. If we knew that the distributions of~$\{\xi^{(i)}\}_{i \geq 1}$ had no atoms, the result would then follow. As we do not have this result, we must resort to other methods to prove Theorem~\ref{thm:leader}.

As we alluded to above, to prove the leadership result what we now need is to show that the distance between the cardinalities of specific blocks goes to infinity. With that in mind, our main goal this section will be to prove the following result:
\begin{proposition} \label{prop:asymp_deg}
For any two distinct indices~$i$ and~$j$, we have
\begin{equation*}
\lim_{t \to \infty} |d_t^{(i)}- d_t^{(j)}| \xrightarrow{t \to \infty} \infty \; \text{almost surely}.
\end{equation*}
\end{proposition}
For now, let us properly show how the above proposition allows us to prove our persistent leadership theorem.
\begin{proof}[Proof of Theorem~\ref{thm:leader}]
Define the events
$$A^{(i)} := \left\{\exists t: \frac{d_t^{(i)}}{\phi_t} > \frac{1}{\sqrt{\phi_i}}\right\},\quad i \in \mathbb{N}.$$
We have
\begin{align*}
\P(A^{(i)}) &= \sum_{n \ge i} \P^{(i)}_n(A^{(i)})\cdot \P(\tau^{(i)} = n) \le \sum_{n \ge i} \exp\left\{-c_1 \cdot \frac{\phi_n}{\sqrt{\phi_i}} \right\}\cdot \P(\tau^{(i)} = n) \le \exp\left\{-c_1 \sqrt{\phi_i}\right\},
\end{align*}
where the first inequality follows from Lemma~\ref{l:upperbounddeg}. Recalling that~$\phi_t \sim b t^\beta$, we then have~$\sum_i \P(A^{(i)}) < \infty$, so by the Borel-Cantelli lemma,
$$\mathbb{P}(A^{(i)} \text{ occurs infinitely often}) = 0.$$
This implies that
$$\exists i^*: \;\frac{d_t^{(j)}}{\phi_t} \le \frac{1}{\sqrt{\phi_j}} \quad \text{for all } j \ge i^* \text{ and }t \in \mathbb{N}.$$
Together with the facts that
$$\frac{1}{\sqrt{\phi_j}} \xrightarrow{j \to \infty} 0 \quad  \text{and}\quad \lim_{t \to \infty} \frac{d_t^{(1)}}{\phi_t} = \zeta_1 > 0,$$ we then obtain that almost surely, there exists a (random) index~$i^{**}$ such that
$$\max\left\{\frac{d_t^{(j)}}{\phi_t}: j \le i^{**}\right\} = \max\left\{\frac{d_t^{(j)}}{\phi_t}: j \in \mathbb{N}\right\}\quad \text{for all }t \text{ large enough},$$
so
$$\max\left\{{d_t^{(j)}}: j \le i^{**}\right\} = \max\left\{{d_t^{(j)}}: j \in \mathbb{N}\right\}\quad \text{for all }t \text{ large enough}.$$
Next, Proposition~\ref{prop:asymp_deg} implies that  for all pairs of distinct~$j,j' \le i^{**}$ we have $$|d_t^{(j)} - d_t^{(j')}| \xrightarrow{t \to \infty} \infty\quad \text{ almost surely.}$$ This implies that there exists~$j^* \le i^{**}$ such that
$$d_t^{(j^*)} > \max\{d_t^{(j)}: j \le i^{**},\;j\neq j^*\} \quad \text{for all } t \text{ large enough}.$$
This concludes the proof.
\end{proof}

\subsection{Increments involving a pair of indices}
For the rest of this section, we fix two indices~$i<j$. In this subsection, we give some definitions and results that will be needed in proving Proposition~\ref{prop:asymp_deg}. Similarly to what was done in Section~\ref{ss:sum_sq_incr}, we will consider the process only at the times~$t$ at which either~$d_t^{(i)}$ or~$d_t^{(j)}$ increases from its previous value. The idea is that, in order to study the distance between cardinalities of blocks~$i$ and~$j$, it is only necessary to look at the process at these times.
\begin{definition}
We define~$\sigma^{(i,j)}_0 := \tau^{(j)} > \tau^{(i)}$ and, for each~$k \in \mathbb{N}_0$,
	$$\sigma^{(i,j)}_{k+1}:= \inf\left\{t > \sigma_k^{(i,j)}: d^{(i)}_t+d^{(j)}_t > d^{(i)}_{\sigma_k^{(i,j)}}+d^{(j)}_{\sigma_k^{(i,j)}}\right\}.$$
	We also let
$$\mathcal{G}^{(i,j)}_k:= \mathcal{F}_{\sigma^{(i,j)}_k},\qquad {D}^{(i,j),i}_k:= d^{(i)}_{\sigma^{(i,j)}_k},\qquad {D}^{(i,j),j}_k:= d^{(j)}_{\sigma^{(i,j)}_k}$$
and
$${D}^{(i,j)}_k:= {D}^{(i,j),i}_k + {D}^{(i,j),j}_k, \qquad Z_k := |D^{(i,j),i}_k - D^{(i,j),j}_k|.$$
\end{definition}
We write~$\Delta D^{(i,j)}_k:= D^{(i,j)}_{k+1}- D^{(i,j)}_k$. 
Recalling the definition of~$\hat{\mathscr{P}}_{m,t}$ in Definition~\ref{def:P_and_P_hat}, as in~\eqref{eq:key_eq_hat_P} we have
\begin{equation*}
\P(\Delta D^{(i,j)}_k = r \mid \mathcal{G}_k^{(i,j)}) = \hat{\mathscr{P}}_{D_k^{(i,j)},\sigma_k^{(i,j)}}(r),\quad r \in \{1,\ldots, R\}.
\end{equation*}
Our goal is to study~$Z_k^{(i,j)}$. To do so, we will first need to understand the process~$D_k^{(i,j)}$, which will control the total amount by which the cardinality of blocks~$i$ and~$j$ may increase at time~$\sigma_k$. Arguing in the same way as in Corollary~\ref{cor:conv_of_scrP}, we obtain:
\begin{lemma}\label{lem:increment_scr_ij}
We have that, almost surely, for each~$r \in \{1,\ldots, R\}$,
$$\P(\Delta D^{(i,j)}_k = r \mid \mathcal{G}^{(i,j)}_k) \xrightarrow{k \to \infty} \begin{cases} \mathscr{P}_0(r) & \text{if } \beta < 1;\\[.2cm] \mathscr{P}_{(\zeta^{(i)}+\zeta^{(j)})b/R}(r)&\text{if }\beta = 1.\end{cases} $$
\end{lemma}

Next, it will be useful to note that, almost surely, for every~$k$ and~$r$,
\begin{equation} \label{eq:bin_r_and_u}
\P(\Delta D^{(i,j),i}_k=r \mid \mathcal{G}_k^{(i,j)}) = \sum_{u= 0}^R \P(\Delta D^{(i,j)}_k = u \mid \mathcal{G}_k^{(i,j)})\cdot \mathscr{B}\left(u,\frac{D^{(i,j),i}_k}{D^{(i,j)}_k},r\right).
\end{equation}
The next two lemmas provide lower bounds for the conditional drift of~$Z^{(i,j)}_k$.
\begin{lemma}
For any~$k$, we almost surely have
\begin{equation}\label{eq:key_ineq_exp}
\E[\Delta Z^{(i,j)}_k \mid \mathcal{G}_k^{(i,j)}]  \ge \frac{Z_k^{(i,j)}}{D_k^{(i,j)}} \cdot \E[\Delta D_k^{(i,j)}\mid \mathcal{G}_k^{(i,j)}]
\end{equation}\end{lemma}
\begin{proof}
On the event~$\{D^{(i,j),i}_k \ge D^{(i,j),j}_k\}$, we have~$Z_k^{(i,j)} = D^{(i,j),i}_k- D^{(i,j),j}_k$, so
\begin{align*}
Z_{k+1}^{(i,j)} &= |D_k^{(i,j),i} + \Delta D_k^{(i,j)i} - D_k^{(i,j),j} - \Delta D_k^{(i,j),j}| \ge Z_k^{(i,j)} + \Delta D^{(i,j),i}_k -\Delta D^{(i,j),j}_k 
\end{align*}
and this gives 
\begin{equation}\label{eq:inter_expec}
\E[\Delta Z^{(i,j)}_k \mid \mathcal{G}_k^{(i,j)}]  \ge \E[\Delta D^{(i,j),i}_k\mid \mathcal{G}_k^{(i,j)}] - \E[\Delta D^{(i,j),j}_k \mid \mathcal{G}_k^{(i,j)}]. 
\end{equation}
Now, using~\eqref{eq:bin_r_and_u} we have
$$\E[\Delta D^{(i,j),i}_k\mid \mathcal{G}_k^{(i,j)}] = \frac{D^{(i,j),i}_k}{D^{(i,j)}_k}\cdot \E[\Delta D_k^{(i,j)}\mid \mathcal{G}_k^{(i,j)}]$$
and similarly for~$j$; using this in~\eqref{eq:inter_expec} we obtain the desired inequality on the event~$\{D^{(i,j),i}_k \ge D^{(i,j),j}_k\}$. Arguing in the same way, we also obtain it in the complementary event~$\{ D_k^{(i,j),j} > D_k^{(i,j),i}\}$.
\end{proof}

\begin{lemma}\label{lem:for_menshikov}
For any~$k$, we almost surely have
\begin{equation}
\E[\Delta Z^{(i,j)}_k \mid \mathcal{G}_k^{(i,j)}] \ge 0
\end{equation}
and
\begin{equation}
\P(\Delta Z^{(i,j)}_k = 0 \mid \mathcal{G}_k^{(i,j)}) < \frac12.
\end{equation}
\end{lemma}
\begin{proof}
The first inequality is an immediate consequence of~\eqref{eq:key_ineq_exp}. The second inequality follows from noting that~$\Delta Z^{(i,j)}_k = 0$ if  and only if~$\Delta D^{(i,j),i}_k = \Delta D^{(i,j),j}_k = \frac12\cdot \Delta D^{(i,j)}_k$, and then using~\eqref{eq:bin_r_and_u} and the elementary observation that~$\mathscr{B}(n,p,n/2) < 1/2$ for any~$n$ and~$p$.
\end{proof}
The next result is the main tool in the application of the Lyapunov function method for the process~$Z_k^{(i,j)}$. It gives a lower bound for the conditional increment of~$Z^{(i,j)}_k$ in terms of~$Z_k^{(i,j)}$ itself.
\begin{proposition}\label{prop:three_quarters}
	Almost surely, there exists~$K_0 \in \N$ such that
	\begin{equation*}
		\E[\Delta Z^{(i,j)}_k \mid \mathcal{G}_k^{(i,j)}] \ge \frac34 \cdot \frac{Z^{(i,j)}_k}{k} \qquad \text{for all }k \ge K_0.
	\end{equation*}
\end{proposition}
\begin{proof}
Let
\begin{equation*}
\mathscr{N}^{(i,j)}:=\begin{cases}1 &\text{if }\beta < 1;\\\sum_{r=1}^R r\cdot \mathscr{P}_{(\zeta^{(i)} + \zeta^{(j)})b/R}(r)&\text{if }\beta = 1.\end{cases}
\end{equation*}
By Lemma~\ref{lem:increment_scr_ij}, we have
\begin{equation*}
\E[\Delta D^{(i,j)}_k \mid \mathcal{G}^{(i,j)}_k] \xrightarrow{k \to \infty} \mathscr{N}^{(i,j)}.
\end{equation*}
Using this and a simple application of the Azuma-Hoeffding inequality, we also obtain
\begin{equation*}
\frac{D^{(i,j)}_k}{k} \xrightarrow{k \to \infty} \mathscr{N}^{(i,j)}.
\end{equation*}
The desired result now follows from combining the two above convergences with~\eqref{eq:key_ineq_exp}.
\end{proof}

\subsection{Transience of difference process}
We are now equipped to take the concluding steps in proving Proposition~\ref{prop:asymp_deg}. The first step in this direction is the following result, which shows that the limit superior of~$Z_k^{(i,j)}$ grows faster than the square root of~$k$.
\begin{proposition} \label{prop:law_of_it_log}
We almost surely have
$$\limsup_{k \to \infty} \frac{Z_k^{(i,j)}}{\sqrt{k}} = \infty. $$
\end{proposition}
\begin{proof}
This is an immediate consequence of Lemma~2 in~\cite{menshikov2008urn}. The assumptions required by that lemma are readily checked using the fact that~$|\Delta Z_k^{(i,j)}| \le R$ and the two inequalities in Lemma~\ref{lem:for_menshikov} (in the notation of~\cite{menshikov2008urn}, take the constant~$a = 1$).  
\end{proof}
We can then bootstrap the above result and finish this section.
\begin{proof}[Proof of Proposition~\ref{prop:asymp_deg}]
It suffices to prove that~$\lim_{k \to \infty}Z_k^{(i,j)} = \infty$ almost surely. To make the notation a bit cleaner, in this proof we will omit the superscript~$(i,j)$, writing~$Z_k$ and~$\mathcal{G}_k$ instead of~$Z_k^{(i,j)}$ and~$\mathcal{G}_k^{(i,j)}$.

For each~$\bar{k} \in \N$, define the event
$$E(\bar{k}):= \left\{\E[\Delta Z_k \mid \mathcal{G}_k] > \frac34\cdot \frac{Z_k}{k} \text{ for all }k \ge \bar{k}\right\}.$$
Note that Proposition~\ref{prop:three_quarters} implies that~$\P(E(\bar{k}))\xrightarrow{\bar{k}\to \infty} 1$. Hence, the proof will be complete if we show that, for any~$\bar{k}$ and any~$\varepsilon > 0$, we have
\begin{equation}\label{eq:want_with_complement}
\mathbb{P}\left(E(\bar{k}) \cap \left\{Z_k \to \infty \right\}^c \right) < \varepsilon. 
\end{equation}

Fix~$\bar{k}$ and~$\varepsilon$; let us prove~\eqref{eq:want_with_complement}. We closely follow the proof of Theorem~1 from~\cite{menshikov2008urn}. Since~$(1+x)^{-2} = 1-2x+3x^2 + o(x^2)$ for~$x$ small, we can choose~$\eta > 0$ such that
\begin{equation} \label{eq:taylor_eta}
(1+x)^{-2} \le 1-2x + 4x^2\quad \text{ for all }x \in [-\eta, \eta].
\end{equation}
Next, define
$$\tau:= \inf\left\{k \ge \bar{k} \vee \frac{R^2}{\eta^2}:\; \frac{k}{Z_k^2} \le \frac{\varepsilon }{16R^2}\right\}.$$
This is a stopping time with respect to~$(\mathcal{G}_{k})_{k \ge 0}$, and by Proposition~\ref{prop:law_of_it_log}, it is almost surely finite. We then let
$$W_\ell := \frac{\tau + \ell}{(Z_{\tau + \ell})^2},\quad \ell \ge 0,$$
with~$W_\ell = \infty$ in case~$Z_{\tau + \ell} = 0$. By the definition of~$\tau$, we have~$W_0 \le \varepsilon /(16R^2)$.

Next, define
\begin{align*}
&L':= \inf\left\{\ell \ge 0:\;W_\ell \ge \frac{1}{16R^2}\right\},\\[.2cm]
&L'':= \inf\left\{\ell \ge 0:\; \mathbb{E}\left[\Delta Z_{\tau + \ell}\mid \mathcal{G}_{{\tau+\ell}}\right] \le \frac34\cdot \frac{Z_{\tau + \ell}}{\tau + \ell}\right\},
\end{align*}
and let~$L:= L' \wedge L''$. These are stopping times with respect to the filtration~$(\mathcal{G}_{{\tau+\ell}})_{\ell \ge 0}$. We now claim that
\begin{equation} \label{eq:claim_super}(W_{\ell \wedge L})_{\ell \ge 0} \text{ is a supermartingale with respect to }(\mathcal{G}_{{\tau+\ell}})_{\ell \ge 0}.
\end{equation} To prove this, we need to check that, for any~$\ell$,
\begin{equation}\label{eq:want_of_Ys}
\mathds{1}\{L > \ell\} \cdot \mathbb{E}\left[\Delta W_{\ell} \mid \mathcal{G}_{\tau+\ell}\right] \le 0.
\end{equation}
To this end, we first write
\begin{equation}\label{eq:for_Delta_E} \begin{split}\Delta W_\ell &= \frac{\tau+\ell +1}{(Z_{\tau+\ell+1})^2} - \frac{\tau+\ell }{(Z_{\tau+\ell})^2}= \frac{\tau + \ell}{(Z_{\tau+\ell})^2}\cdot \left[\left(1+\frac{1}{\tau+\ell}\right)\cdot \left(\frac{1}{1+\frac{\Delta Z_{\tau+\ell}}{Z_{\tau+\ell}}} \right)^2 - 1\right].\end{split}  \end{equation}
Now, we have~$|\Delta Z_{\tau + \ell}| \le R$ and on the event~$\{L > \ell\}$ we have, by the definition of~$L'$,
$$W_\ell = \frac{\tau+\ell}{(Z_{\tau+\ell})^2} < \frac{1}{16R^2} < 1, \quad \text{so}\quad Z_{\tau+\ell} > (\tau+\ell)^{1/2} > \frac{R}{\eta},$$
where in the last inequality we have used that~$\tau > \frac{R^2}{\eta^2}$. This shows that on the event~$\{L > \ell\}$ we have~$\frac{\Delta Z_{\tau + \ell}}{Z_{\tau + \ell}} \le \eta$, so, using~\eqref{eq:taylor_eta},
$$\left(\frac{1}{1+\frac{\Delta Z_{\tau+\ell}}{Z_{\tau+\ell}}} \right)^2\; {\le} 1 - 2\cdot \frac{\Delta Z_{\tau + \ell}}{Z_{\tau + \ell}} + 4\cdot \left(\frac{\Delta Z_{\tau + \ell}}{Z_{\tau + \ell}}\right)^2 \le 1 - 2\cdot \frac{\Delta Z_{\tau + \ell}}{Z_{\tau + \ell}} +  \frac{4 R^2}{(Z_{\tau + \ell})^2}.$$
Plugging this into~\eqref{eq:for_Delta_E}, we obtain that the left-hand side of~\eqref{eq:want_of_Ys} is bounded from above by
\begin{equation}\label{eq:exp_sqbk}
\mathds{1}\{L > \ell\} \cdot  \frac{\tau + \ell}{(Z_{\tau+\ell})^2}\cdot \left[\left(1+\frac{1}{\tau+\ell}\right)\cdot \left( 1 - 2\cdot \frac{\mathbb{E}[\Delta Z_{\tau + \ell}\mid \mathcal{G}_{{\tau + \ell}}]}{Z_{\tau+\ell}} + \frac{4R^2}{(Z_{\tau +\ell})^2}\right) - 1\right].\end{equation}
Now, on~$\{L > \ell\}$ we have~$\mathbb{E}[\Delta Z_{\tau+\ell} \mid \mathcal{G}_{{\tau+\ell}}] > \frac34\cdot \frac{Z_{\tau+\ell}}{\tau + \ell}$ (by the definition of~$L''$) and
$$W_\ell = \frac{\tau+\ell}{(Z_{\tau+\ell})^2} < \frac{1}{16R^2}, \quad \text{so}\quad \frac{4R^2}{(Z_{\tau+\ell})^2} < \frac{1}{4(\tau + \ell)}$$
(by the definition of~$L'$).
  This shows that the expression inside the square brackets in~\eqref{eq:exp_sqbk} is bounded from above by
\begin{align*}&\left(1+\frac{1}{\tau+\ell}\right)\cdot \left( 1 -  \frac{3}{2(\tau+\ell)} + \frac{1}{4(\tau+\ell)}\right) - 1= \left(1+\frac{1}{\tau+\ell}\right)\cdot \left(1- \frac{5}{4(\tau+\ell)}\right) -1 < 0. \end{align*}
This concludes the proof of~\eqref{eq:claim_super}.

Since it is a non-negative supermartingale,~$(W_{\ell \wedge L})$ converges almost surely to a non-negative random variable~$W_\infty$. We have~$\mathbb{E}[W_\infty] \le \mathbb{E}[W_0] < \frac{\varepsilon }{16R^2}$, so~$W_\infty$ is almost surely finite. Moreover, we have
\begin{align*}\frac{\varepsilon}{16R^2} \ge \mathbb{E}[W_0] \ge \mathbb{E}[W_L] &\ge \mathbb{E}[W_L \cdot \mathds{1}\{L'<\infty,\; L' < L''\}]\ge \frac{1}{16R^2}\cdot \mathbb{P}(L'< \infty,\; L' < L''),\end{align*}
so
\begin{equation}\label{eq:on_lp}\mathbb{P}(L'< \infty,\; L' < L'') < \varepsilon.\end{equation}
Next, note that 
\begin{equation}\{L'' < \infty\} \subseteq E(\bar{k})^c, \label{eq:on_lpp}\end{equation}
by the definition of~$L''$ and of~$E(\bar{k})$. Finally, note that on~$\{L = \infty\}$, we have~$W_{\ell \wedge L} = W_\ell = \frac{\tau+\ell}{(Z_{\tau + \ell})^2}$; the fact that this converges to a finite limit on~$\{L = \infty\}$ then implies that
\begin{equation} \label{eq:on_conv_of_dk}Z_k \to \infty \text{ on }\{L = \infty\}.\end{equation}
Putting~\eqref{eq:on_lp},~\eqref{eq:on_lpp} and~\eqref{eq:on_conv_of_dk} together now gives the desired inequality~\eqref{eq:want_with_complement}.
\end{proof}

\section{Convergence and central limit theorem for the maximum}\label{sec:maximum}
%!TEX root = ./ms.tex
We can now finish the first part of our results regarding the maximum cardinality. Intuitively, the main idea is that, for large time~$t$, the maximum cardinality \emph{is} the cardinality of a given block, so that Theorem~\ref{thm:convergence} applies.
\begin{proof}[Proof of Theorem~\ref{thm:max}, convergence]
    With Theorem~\ref{thm:leader} at hand, it is quite easy to prove that
    \begin{equation}\label{eq:easy_conv}\lim_{\substack{t \to \infty\\\mathrm{a.s.}}}\frac{1}{t^\beta}\cdot \max_{i}d_t^{(i)} = \lim_{\substack{t \to \infty\\\mathrm{a.s.}}}\frac{1}{t^\beta}\cdot d_t^{(I)} = \xi^{(I)} = \sup_i \xi^{(i)} \in (0,\infty). \end{equation}
    Indeed, for two distinct indices~$i,j$, on the event~$\{\xi^{(i)} > \xi^{(j)}\}$ we have~$d_t^{(i)} - d_t^{(j)} \to \infty$ (since~$d_t^{(i)}/t^\beta \to \xi^{(i)}$ and~$d_t^{(j)}/t^\beta\to \xi^{(j)}$). Hence, the fact that~$d_t^{(I)} - \max_{i \neq I} d_t^{(i)} \to \infty$ gives~\eqref{eq:easy_conv}.
    
    To prove the convergence in~$L_p$ for all~$p \in [1,\infty)$, let
    $$M_t := \max_{i \ge 1} \frac{d_t^{(i)}}{\phi_t},\;t \ge 1\qquad \text{and} \qquad M_\infty := \lim_{\substack{t\to \infty\\\mathrm{a.s.}}}M_t = \sup_i \zeta^{(i)}.$$
    Using the facts, given in Lemma~\ref{lem:x_mart} and Proposition~\ref{prop:zetaasconv}, that~$\{d_t^{(i)}/\phi_t\}_{t \ge 1}$ is a submartingale with $d_t^{(i)}/\phi_t \xrightarrow{t \to \infty} \zeta^{(i)}$ almost surely (for any~$i$), we obtain the bound, for any~$p \in \mathbb{N}$, with \new{$p > 1/\beta$}:
    \begin{align}
    \label{eq:moment_super_1}& \mathbb{E}[(M_t)^p] \le \sum_{i \ge 1}\mathbb{E}[(d_t^{(i)}/\phi_t)^p] \le \sum_{i \ge 1}\mathbb{E}[(\zeta^{(i)})^p] \stackrel{\eqref{eq:zetakmoment}}{<} \infty.
    \end{align}
    Again using the fact  that~$\{d_t^{(i)}/\phi(t)\}_{t \ge 1}$ is a submartingale, it is easy to check that $\{M_t\}_{t \ge 1}$ is also a submartingale. The above bound tell us that it is bounded in~$L_p$ for every~$p \in \mathbb{N}$ (hence also for every~$p \in [1,\infty)$). It then follows that~$M_t \to M_\infty$ in~$L_p$, for~$p \in [1,\infty)$.
    \end{proof}
    
    Before proving the remaining part of Theorem~\ref{thm:max}, we give a definition and state and prove two lemmas. The definition concerns the index of a block with maximum cardinality at time~$n$. The lemmas show respectively that a CLT is valid for this block that \emph{had} maximum cardinality at time~$n$, and that the law of the fluctuations of the maximum cardinality at time~$n$ converges weakly, as~$n$ goes to infinity, to the law~$\mu^*$ defined in Theorem~\ref{thm:max}.
    \begin{definition}
    For each~$n \in \N$, let
    $$I_n:= \min\{i: d^{(i)}_n = \max_j d^{(j)}_n\}.$$
    \end{definition}
    Recall the definition of the random variable~$I$ that appears in Theorem~\ref{thm:leader}. The statement of that theorem implies that
    \begin{equation}\label{eq:conv_in_i}
    \P(I_n \neq I) \xrightarrow{n \to \infty} 0.
    \end{equation}
    \begin{lemma}\label{lem:first_conv_mu}
    We have that
    \begin{equation}
    t^{\beta/2} \cdot \left(\frac{d_t^{(I_n)}}{t^\beta} - \xi^{(I_n)} \right) \xrightarrow[\mathrm{(d)}]{t \to \infty} \mu^*_n,
    \end{equation}
    where~$\mu^*_n$ is the distribution of~$W \cdot Z^*_n$, where~$W, Z^*_n$ are independent,~$W$ is a standard Gaussian and
    $$(Z^*_n)^2 \stackrel{\mathrm{(d)}}{=} \begin{cases}
    \xi^{(I_n)}&\text{if }\beta<1;\\ \xi^{(I_n)}\cdot \left(1-\frac{\xi^{(I_n)}}{R}\right)&\text{if }\beta = 1.
    \end{cases} $$
    \end{lemma}
    \begin{proof}
    Define, for each~$i,n \in \N$,
    $$\tilde{\P}^{(i)}_n(\cdot):= \P(\cdot \mid I_n = i),$$
    and~$\tilde{\E}^{(i)}_n$ be the associated expectation operator.
    Recall that~$\frac{d_t^{(i)}}{\phi_t} \to \zeta^{(i)}$ almost surely, so this convergence also holds~$\tilde{\P}^{(i)}_n$-almost surely. Moreover, repeating the proof of Proposition~\ref{prop:e_incr_d_sq} shows that 
    $$\phi_t \cdot \sum_{s \ge t} \tilde{\E}^{(i)}_n[(\Delta X_s^{(i)})^2] \xrightarrow{t \to \infty} \begin{cases}
    \tilde{\E}^{(i)}_n[\zeta^{(i)}]& \text{if }\beta<1;\\[.2cm]
    \tilde{\E}^{(i)}_n\left[\zeta^{(i)}\cdot \left( 1- \frac{\zeta^{(i)}\new{b}}{R}\right) \right]&\text{if }\beta=1.
    \end{cases}$$
    Putting these convergences together as in the proof of Theorem~\ref{thm:clt}, we obtain that, under~$\tilde{\P}^{(i)}_n$,
    $$t^{\beta/2}\cdot \left(\frac{d_t^{(i)}}{t^\beta} - \xi^{(i)}\right) \xrightarrow[\mathrm{(d)}]{t \to \infty} \tilde{\mu}^{(i)}_n,$$
    where~$\tilde{\mu}^{(i)}_n$ is the distribution of~$W\cdot \tilde{\mathcal{Z}}^{(i)}_n$, where~$W,\tilde{\mathcal{Z}}^{(i)}_n$ are independent,~$W$ is a standard Gaussian and
    $$(\tilde{\mathcal{Z}}^{(i)}_n)^2 \sim \begin{cases}
    \text{law of }\xi^{(i)} \text{ under }\tilde{\P}^{(i)}_n&\text{if }\beta < 1;\\[.2cm]
    \text{law of }\xi^{(i)}\cdot \left(1-\frac{\xi^{(i)}}{R}\right) \text{ under }\tilde{\P}^{(i)}_n&\text{if }\beta = 1.
    \end{cases} $$
    Now, as in the proof of Theorem~\ref{thm:clt}, the proof is concluded by observing that~$\mu^{*}_n = \sum_{i} \P(I_n = i)\cdot \tilde{\mu}^{(i)}_n$.
    \end{proof}
    Recall the definition of the distribution~$\mu^*$ in Theorem~\ref{thm:max}, and the definition of~$\mu_n^*$ in Lemma~\ref{lem:first_conv_mu}. The following result is an easy consequence of~\eqref{eq:conv_in_i}; the proof is omitted. 
    \begin{lemma}\label{lem:mu_n_star}
    As~$n \to \infty$,~$\mu^*_n$ converges weakly to~$\mu^*$.
    \end{lemma}
    We can now finish the proof of the result.
    \begin{proof}[Proof of Theorem~\ref{thm:max}, central limit theorem]
    Let~$h: \R \to \R$ be a continuous and bounded function. Given~$\varepsilon > 0$, using~\eqref{eq:conv_in_i} and Lemma~\ref{lem:mu_n_star} choose~$n$ large enough that
    $$\P(I_n \neq I) < \frac{\varepsilon}{3\|h\|_\infty}\qquad \text{and}\qquad \left|\int h\;\mathrm{d}\mu^*_n - \int h\; \mathrm{d}\mu^*\right| < \frac{\varepsilon}{3}. $$
    Next, using Lemma~\ref{lem:first_conv_mu}, choose~$t$ large enough that
    $$\left| \E\left[ h\left( t^{\beta/2} \cdot \left(\frac{d_t^{(I_n)}}{t^\beta} - \xi^{(I_n)} \right)  \right)\right] - \int h\;\mathrm{d}\mu_n^* \right|<\frac{\varepsilon}{3}. $$
    Finally note that
    $$\left| \E\left[ h\left( t^{\beta/2} \cdot \left(\frac{d_t^{(I_n)}}{t^\beta} - \xi^{(I_n)} \right)  \right)\right] - \E\left[ h\left( t^{\beta/2} \cdot \left(\frac{d_t^{(I)}}{t^\beta} - \xi^{(n)} \right)  \right)\right] \right|  < \|h\|_\infty \cdot \P(I_n \neq I) < \frac{\varepsilon}{3}.  $$
    Putting things together, we conclude that
    $$\E\left[ h\left( t^{\beta/2} \cdot \left(\frac{d_t^{(I)}}{t^\beta} - \xi^{(n)} \right)  \right)\right] \xrightarrow{t \to \infty} \int h\;\mathrm{d}\mu^*,$$
    as required.
    \end{proof}

{\bf Acknowledgements }  \textit{C.A.}  was partially supported by  the Deutsche Forschungsgemeinschaft (DFG), and the Noise-Sensitivity everywhere ERC Consolidator Grant 772466.  \textit{R.R.} was supported by the project Stochastic Models of Disordered and Complex Systems. The Stochastic Models of Disordered and Complex Systems is a Millennium Nucleus (NC120062) supported by the Millenium Scientific Initiative of the Ministry of Science and Technology  (Chile).

\appendix
\section{Appendix}
%!TEX root = ./ms.tex
\subsection{Martingale concentration inequalities}
For the sake of completeness, we state here two useful concentration inequalities for martingales which are used throughout the paper.
\begin{theorem}[Azuma-Hoeffding Inequality -- \cite{CLBook}]\label{thm:azuma} Let $(M_n,\mathcal{F}_n)_{n \ge 1}$ be a martingale. Assume there exists a sequence of negative real numbers~$(a_n)_{n \ge 1}$ such that~$|M_{n+1} - M_n| \le a_n$ for each~$n$.
 Then, 
	\[
	\P \left( | M_n - M_0 | > \lambda \right) \le \exp\left\{ -\frac{\lambda^2}{\sum_{i=1}^n a_i^2} \right\} \quad \text{for all }\lambda > 0,\;n\in\mathbb{N}.
	\]
\end{theorem}

\begin{theorem}[Freedman's Inequality -- \cite{F75}]\label{thm:freedman} Let $(M_n, \mathcal{F}_n)_{n \ge 1}$ be a (super)martingale. Assume that~$M_0 = 0$ and there exists~$K > 0$ such that~$|M_{n+1}-M_n| \le K$ for all~$n$.
	Write 
	\[
	V_n := \sum_{k=1}^{n-1} \mathbb{E} \left[(M_{k+1}-M_k)^2 \mid\mathcal{F}_k \right],\quad n \in \N.
	\]
	Then, 
	\[
	\P \left(\exists n:\;   M_n \ge \lambda\text{ and } V_n \le \sigma^2\right) \le \exp\left\{ -\frac{\lambda^2}{2\sigma^2 + 2K\lambda /3} \right\}\quad \text{for all }\lambda > 0.
	\]
\end{theorem}

\subsection{Martingale central limit theorem} The CLT results present in this paper follow from an application of the following CLT concerning the tails of converging martingales.
\label{ss:clt_martingale}
\begin{theorem}[Martingale central limit theorem -- \cite{MartingaleCLTbook2}, Corollary~3.5, page 79] \label{thm:clt_martingale}
Let~$(S_n)_{n \in \mathbb{N}}$ be a square-integrable martingale with respect to a filtration~$(\mathscr{F}_n)_{n \in \N}$, satisfying
\begin{equation}\label{eq:mart_clt_1}
\sum_{n=1}^\infty \E[(\Delta S_n)^2] < \infty.
\end{equation}	
Let~$s_n := \left(\sum_{m=n}^\infty \mathbb{E}[(\Delta S_m)^2]\right)^{1/2} $. Assume that
\begin{equation}
\label{eq:mart_clt_2}
\frac{1}{s_n}\cdot \sup_{m \ge n}|\Delta S_m| \xrightarrow[\mathrm{prob.}]{n \to \infty} 0,
\end{equation}
that
\begin{equation} \label{eq:mart_clt_22}
\frac{1}{(s_n)^2} \cdot\E\left[  \sup_{m \ge n}(\Delta S_m)^2\right] < \infty,
\end{equation}
and that
\begin{equation} \label{eq:mart_clt_222}\frac{1}{(s_n)^2}\cdot \sum_{m=n}^\infty (\Delta S_m)^2 \xrightarrow[\mathrm{prob.}]{n \to \infty} \eta^2\end{equation} for some random variable~$\eta^2$. Then,~$\frac{1}{s_n}\cdot \sum_{m=n}^\infty \Delta S_m$ converges in distribution to the probability distribution with characteristic function equal to~$t \mapsto \E\left[\exp\{-\eta^2 t^2/2\}\right]$.
\end{theorem}

\subsection{Proofs of results on series along sets with asymptotic density} We prove here the lemmas about sums along sets with asymptotic density used to establish the CLT results present in this paper.
\begin{proof}[Proof of Lemma~\ref{lem:sum_over}]
The result will readily follow once we show that, for every integer~$K \ge 2$, we have
\begin{equation}\label{eq:want_forgnnu}
\lim_{N \to \infty} N \sum_{n \in \Lambda \cap [N,KN]}\frac{1}{n^2} = \alpha\left(1-\frac{1}{K}\right).
\end{equation}
To that end, fix~$K$ and define, for each~$N \in \mathbb{N}$, the function
$$g(u):= \lfloor Nu \rfloor^{-2},\quad u \in [1,K]$$
and the measure
$$\nu_N^\Lambda:= \sum_{n \in \Lambda \cap [N,KN]} \delta_{\{n/N\}} $$
on Borel sets of~$[1,K]$. Note that
\begin{equation}\label{eq:funny_sum}
\sum_{n \in \Lambda \cap [N,KN]}\frac{1}{n^2} = \int_{[1,K]} g_N\;\mathrm{d}\nu^\Lambda_N.
\end{equation}
It is also readily seen that
\begin{equation}\label{eq:conv_gn}N^2\cdot g_N(u) \xrightarrow{N \to \infty} \frac{1}{u^2} \text{ uniformly on }[1,K].\end{equation}
Using the fact that~$\Lambda$ has asymptotic density~$\alpha$, it is a routine exercise to show that
\begin{equation} \label{eq:conv_nu} \frac{1}{N}\cdot\nu_N^\Lambda([1,u]) \xrightarrow{N \to \infty} \alpha (u-1) \text{ uniformly on }[1,K].\end{equation}

Now, writing~$h(u):=u^{-2}$ and letting~$\ell$ denote the Lebesgue measure on~$[1,K]$, we bound
\begin{align*}
&\left|N\sum_{n \in \Lambda \cap [N,KN]} \frac{1}{n^2} - \alpha \left(1 - \frac{1}{K} \right) \right| \stackrel{\eqref{eq:funny_sum}}{=} \left| \int_{[1,K]} (N^2g_N)\; \mathrm{d}\left(\frac{1}{N}\nu^\Lambda_N\right)  - \int_{[1,K]} h\; \mathrm{d}(\alpha\ell)\right|\\[.2cm]
&\le \int_{[1,K]} (N^2 \cdot g_N)\;\mathrm{d}\left|\frac{1}{N}\nu^\Lambda_N - \alpha \ell\right|   + \int_{[1,K]} \left|N^2\cdot g_N - h \right|\mathrm{d}(\alpha \ell) \xrightarrow{N \to \infty} 0,
\end{align*}
the convergence following from~\eqref{eq:conv_gn} and~\eqref{eq:conv_nu}. This completes the proof of~\eqref{eq:want_forgnnu}.
\end{proof}

\begin{proof}[Proof of Lemma~\ref{lem:sop_ad}]
Let~$m:= \sum_r r\rho_r$. It is easy to check that the assumptions imply that
\begin{equation}\label{eq:conv_ak}
\lim_{k \to \infty} \frac{a_k}{k} = m.
\end{equation}
Define, for~$r \in \{1,\ldots, R\}$, the sets
$$\Lambda_r:= \{k\in\mathbb{N}:\Delta a_k = r\},\qquad \Xi_r := \{a_k: k \in \Lambda_r\}.$$
We claim that, for each~$r$,~$\Xi_r$ has asymptotic density~$\rho_r/m$. To check this, fix~$\varepsilon > 0$. Then,~\eqref{eq:conv_ak} implies that, for~$N$ large enough,
$$\Xi_r \cap \left\{1,\ldots, a_{\left \lfloor (1-\varepsilon)N/m\right\rfloor}\right\} \subset \Xi_r \cap \{1,\ldots,N\} \subset \Xi_r \cap \left\{1,\ldots, a_{\left \lceil (1+\varepsilon)N/m\right\rceil}\right\},$$
so
$$\left|\Lambda_r \cap \left\{1,\ldots, \left \lfloor \frac{(1-\varepsilon)N}{m} \right \rfloor \right\}\right| \le |\Xi_r \cap \{1,\ldots, N\}| \subset \left|\Lambda_r \cap \left\{1,\ldots, \left \lceil \frac{(1+\varepsilon)N}{m} \right \rceil \right\}\right| $$
and then, again by~\eqref{eq:conv_ak},
$$(1-\varepsilon)\frac{\rho_r}{m} < \liminf_{N \to \infty} \frac{|\Xi_r \cap \{1,\ldots, N\}|}{N} \le  \limsup_{N \to \infty} \frac{|\Xi_r \cap \{1,\ldots, N\}|}{N} \le(1+\varepsilon)\frac{\rho_r}{m},$$
so, since~$\varepsilon$ is arbitrary, the claim is proved.

We now write
$$a_{k_0}\sum_{k \ge k_0} \left(\frac{\Delta a_k}{a_k}\right)^2 = \sum_{r=1}^R r^2\cdot a_{k_0}\sum_{\substack{k\ge k_0,\\ k\in \Lambda_r}}\frac{1}{(a_k)^2} = \sum_{r=1}^R  r^2 \cdot a_{k_0} \sum_{\substack{b \ge a_{k_0},\\ b \in \Xi_r}}\frac{1}{b^2} \xrightarrow{k_0 \to \infty} \sum_{r=1}^R r^2 \cdot \frac{\rho_r}{m},$$
where the convergence follows from~\eqref{eq:conv_ak}, the fact that~$\Xi_r$ has asymptotic density~$\rho_r/m$ and Lemma~\ref{lem:sum_over}. This completes the proof.
\end{proof}

\bibliographystyle{plain}
\bibliography{ref}
\end{document}